\documentclass{article}
\usepackage[margin=1.25in]{geometry}
\usepackage{graphicx}
\usepackage{lipsum}
\usepackage{amsfonts}
\usepackage{amsmath}
\usepackage{multirow}
\usepackage{graphicx}
\usepackage{epstopdf}
\usepackage{algorithm}
\usepackage{algorithmic}
\usepackage{marvosym}
\usepackage{bm}
\usepackage{booktabs}
\usepackage{longtable}
\usepackage{amssymb}

\newcommand{\x}{\mathbf{x}}
\newcommand{\y}{\mathbf{y}}
\newcommand{\z}{\mathbf{z}}
\newcommand{\w}{\mathbf{w}}
\newcommand{\e}{\mathbf{e}}

\renewcommand{\d}{\mathbf{d}}

\renewcommand{\v}{\mathbf{v}}
\newcommand{\s}{\mathbf{s}}
\renewcommand{\a}{\mathbf{a}}
\renewcommand{\c}{\mathbf{c}}
\renewcommand{\b}{\mathbf{b}}

\newcommand{\Z}{\mathbf{Z}}
\newcommand{\W}{\mathbf{W}}
\newcommand{\A}{\mathbf{A}}

\renewcommand{\S}{\mathbf{S}}

\newcommand{\I}{\mathbf{I}}
\newcommand{\U}{\mathbf{U}}
\newcommand{\V}{\mathbf{V}}

\renewcommand{\S}{\mathbf{S}}

\newcommand{\<}{\left\langle}
\renewcommand{\>}{\right\rangle}

\newtheorem{theorem}{Theorem}
\newtheorem{proposition}{Proposition}

\newtheorem{lemma}{Lemma}
\newtheorem{proof}{Proof}
\newtheorem{definition}{Definition}

\DeclareMathOperator*{\argmin}{argmin}

\begin{document}

\title{Accelerated Alternating Direction Method of Multipliers: an Optimal $O(1/K)$ Nonergodic Analysis
}


\author{Huan Li\footnote{Peking University. Email: lihuanss@pdu.edu.cn} \and
Zhouchen Lin\footnote{Peking University. Email: zlin@pku.edu.cn}
}

\maketitle

\begin{abstract}
  The Alternating Direction Method of Multipliers (ADMM) is widely used for linearly constrained convex problems. It is proven to have an $o(1/\sqrt{K})$ nonergodic convergence rate and a faster $O(1/K)$ ergodic rate after ergodic averaging, where $K$ is the number of iterations. Such nonergodic convergence rate is not optimal. Moreover, the ergodic averaging may destroy the sparseness and low-rankness in sparse and low-rank learning. In this paper, we modify the accelerated ADMM proposed in [Y. Ouyang, Y. Chen, G. Lan, and E. Pasiliao, An Accelerated Linearized Alternating Direction Method of Multipliers, SIAM J. on Imaging Sciences, 2015, 1588-1623] and give an $O(1/K)$ nonergodic convergence rate analysis, which satisfies $|F(\x^K)-F(\x^*)|\leq O(1/K)$, $\|\A\x^K-\b\|\leq O(1/K)$ and $\x^K$ has a more favorable sparseness and low-rankness than the ergodic peer, where $F(\x)$ is the objective function and $\A\x=\b$ is the linear constraint. As far as we know, this is the first $O(1/K)$ nonergodic convergent ADMM type method for the general linearly constrained convex problems. Moreover, we show that the lower complexity bound of ADMM type methods for the  separable linearly constrained nonsmooth convex problems is $O(1/K)$, which means that our method is optimal.
\end{abstract}


\section{Introduction}

We consider the following general linearly constrained convex problem:
\begin{equation}\label{problem}
\min_{\x_i \in \mathbb{R}^{n_i}} \sum_{i=1}^2 \left(f_i(\x_i)+h_i(\x_i)\right),\quad s.t.  \quad \sum_{i=1}^2\A_i\x_i=\b,
\end{equation}
where both $f_i$ and $h_i$ are convex. $f_i$ is $L_i$-smooth and $h_i$ can be nonsmooth. Specially, $f_i$ can vanish in problem (\ref{problem}). Problems like (\ref{problem}) arise from diverse applications in machine learning, imaging and computer vision, see, e.g., \cite{Boyd-2010-Distributed,Chambolle-2011,Esser-2010} and references therein. In machine learning, $f_i$ is often the loss function to fit the data and $h_i$ is the regularizer that promotes some prior information on the desired solution, such as sparseness and low-rankness. We say $f_i$ is $L_i$-continuous if it satisfies $|f_i(\x_i)-f_i(\y_i)|\leq L_i\|\x_i-\y_i\|,\forall \x_i,\y_i$, and $L_i$-smooth if $\nabla f_i$ is $L_i$-continuous: $\|\nabla f_i(\x_i)-\nabla f_i(\y_i)\|\leq L_i\|\x_i-\y_i\|,\forall \x_i,\y_i$. We denote $F_i(\x_i)=f_i(\x_i)+h_i(\x_i)$, $\x=(\x_1,\x_2)$, $F(\x)=\sum_{i=1}^2 F_i(\x_i)$ and $\A\x=\sum_{i=1}^2\A_i\x_i$. The discussion in this paper also suits for the general constraint $\sum_{i=1}^2\mathcal{A}_i(\x_i)=\b$, where $\mathcal{A}_i:\mathbb{R}^{n_i}\rightarrow \mathbb{R}^{m}$ is a linear mapping. For simplicity we focus on $\sum_{i=1}^2\A_i\x_i=\b$. We denote $\|\x\|$ as $\|\x\|_2$ for a vector $\x$.

ADMM \cite{Boyd-2010-Distributed} is widely used in imaging and vision to solve problem (\ref{problem}) since the separable structure can be exploited. ADMM consists of three steps:
\begin{subequations} \begin{align}
&\x_1^{k+1}=\argmin_{\x_1} L(\x_1,\x_2^k,\lambda^k,\rho),\label{adm1}\\
&\x_2^{k+1}=\argmin_{\x_2} L(\x_1^{k+1},\x_2,\lambda^k,\rho),\label{adm2}\\
&\lambda^{k+1}=\lambda^k+\rho\left(\sum_{i=1}^2\A_i\x_i^{k+1}-\b\right),\label{adm3}
\end{align} \end{subequations}
where
\begin{equation}
L(\x_1,\x_2,\lambda,\rho)=\sum_{i=1}^2 (f_i(\x_i)+h_i(\x_i))+\<\lambda,\sum_{i=1}^2\A_i\x_i-\b\>+\frac{\rho}{2}\left\|\sum_{i=1}^2\A_i\x_i-\b\right\|^2\notag
\end{equation}
is the augmented Lagrangian function and $\lambda$ is the Lagrange multiplier. When $F_i$ is not simple and $\A_i$ is non-unitary, the cost of solving the subproblems may be high. Thus the Linearized ADMM (LADMM) is proposed by linearizing the augmented term $\|\A\x-\b\|^2$ and the complex $f_i$ \cite{he-2002,shefi-2014,Wang-2012-LADM} such that the subproblems may even have closed form solutions.

Traditional convergence rate analysis on ADMM is difficult due to its serial update of $\x_1$ and $\x_2$, which means that $(\x_1^{k+1},\x_2^{k+1})$ is not the solution to $\min_{\x_1,\x_2} L(\x_1,\x_2,\lambda^k)$. Thus some alternative criteria are used instead. The most popular criterion is the ergodic convergence rate.
\begin{definition}
Let $\{\x^1,\cdots,\x^K\}$ be a sequence produced by the algorithm that they have the property promoted by the regularizer $h(\x)$, for instance, sparseness and low-rankness. We say a convergence rate is nonergodic if it measures the optimality at $\x^K$ directly. A convergence rate is ergodic if it considers the optimality at the point of $\sum_{k=1}^K c_k\x^k$ with $c_k> 0$ and $\sum_{k=1}^K c_k=1$.
\end{definition}

The most commonly used ergodic criterion for ADMM is the average form of $\frac{1}{K}\sum_{k=1}^K\x^k$. It is proved in \cite{He-2011-Rate} that ADMM converges with an $O(1/K)$ ergodic rate. A critical disadvantage of the ergodic result is that the point measured for the convergence rate analysis may not have the property promoted by $h(\x)$ since it may be destroyed by the ergodic averaging. For example, in sparse learning, $\{\x^1,\cdots,\x^K\}$ are sparse, but their average may not be sparse any more. So the nonergodic analysis is strongly required for ADMM. He and Yuan \cite{He-nonergodic-ADM} proved $\|\w^{K+1}-\w^K\|^2\leq \frac{1}{K}$ with $\w^K=(\x^K,\lambda^K)$. However, this criterion does not directly measure how far $F(\x^K)$ is from $F(\x^*)$ and how much the constraint error $\|\A\x^K-\b\|$ is, where $\x^*$ is an optimal solution to problem (\ref{problem}).

Recently, Davis and Yin \cite{Davis-2014-DR} proved that the Douglas-Rachford (DR) splitting \cite{Douglas-1956-DR} converges with an $O(1/K)$ ergodic rate and an $o(1/\sqrt{K})$ nonergodic rate. Moreover, they constructed some examples showing that this rate is tight. As is known, ADMM is a special case of DR splitting \cite{Gabay-1983-DR}. So for ADMM, Davis and Yin \cite{Davis-2014-DR} established $|F(\x^K)-F(\x^*)|\leq o(1/\sqrt{K})$ and $\|\A\x^K-\b\|\leq o(1/\sqrt{K})$ in a nonergodic sense. Thus in sparse and low-rank learning, we have that for ADMM the nonergodic solution $\x^K$ is sparse or low-rank, but has the slow $o(1/\sqrt{K})$ theoretical convergence rate, and the ergodic solution $\sum_{k=1}^K c_k\x^k$ has the faster $O(1/K)$ theoretical convergence rate, but may not be sparse or low-rank. We want to combine the advantages of these two aspects, i.e., a faster $O(1/K)$ convergence rate but still in the nonergodic sense. This paper aims to solve this problem via using Nesterov's acceleration scheme for ADMM.

Beck and Teboulle \cite{Beck-2009-APG} extended Nesterov's accelerated gradient method \cite{Nesterov1983} to the nonsmooth unconstrained problem of $\min_{\x} f(\x)+h(\x)$, which consists of two steps: first extrapolates a point $\y^k=\x^k+\frac{\theta^k(1-\theta^{k-1})}{\theta^{k-1}}(\x^k-\x^{k-1})$ and then computes $\x^{k+1}=\mbox{Prox}_{\alpha h}\left(\y^k-\alpha\nabla f(\y^k)\right)$, where $\mbox{Prox}_{\alpha h}(\z)=\argmin_{\x} h(\x)+\frac{1}{2\alpha}\|\x-\z\|^2$. On the other hand, Nesterov \cite{Nesterov1988} proposed another accelerated gradient method, which consists of three steps: $\z^k=(1-\theta^{k-1})\z^{k-1}+\theta^{k-1}\x^k$, $\y^k=(1-\theta^k)\z^k+\theta^k\x^k$ and $\x^{k+1}=\mbox{Prox}_{\frac{\alpha}{\theta^k}h}\left(\x^k-\frac{\alpha }{\theta^k}\nabla f(\y^k)\right)$. We follow \cite{Tseng-2008} to name these two schemes as Nesterov's first and second acceleration scheme, respectively.

Chen et al. \cite{Chen-2015-IPADM} proposed an inertial proximal ADMM which uses the same idea as Nesterov's first scheme: first extrapolates a point $(\hat\x_1^k,\hat\x_2^k,\hat\lambda^k)$ and then performs the steps (\ref{adm1})-(\ref{adm3}) on $(\hat\x_1^k,\hat\x_2^k,\hat\lambda^k)$. However, they only established the $o(1/\sqrt{K})$ convergence rate in the sense of $\min_{k=1,\cdots,K}|F(\x^k)-F(\x^*)|\leq o(1/\sqrt{K})$ and $\min_{k=1,\cdots,K} \|\A\x^k-\b\|\leq o(1/\sqrt{K})$. Lorenz and Pock \cite{Lorenz-2015} analyzed the inertial forward-backward algorithm for the general  monotone inclusions, which include problem (\ref{problem}) as a special case. However, no convergence rate is established in \cite{Lorenz-2015}.

Ouyang et al. \cite{ouyang-2015-fastLADM} proposed an accelerated ADMM via Nesterov's second acceleration scheme. The convergence rate is better than that of LADMM in terms of their dependence on the Lipschitz constant of the smooth component. However, the entire convergence rate remains $O(1/K)$ in an ergodic sense. Nesterov's second scheme only influences the linearization of $f_i$ in steps (\ref{adm1})-(\ref{adm2}). It cannot improve the nonergodic rate of ADMM. Thus, the nonergodic rate of the accelerated ADMM in \cite{ouyang-2015-fastLADM} cannot be better than $o(1/\sqrt{K})$. Please see Section \ref{review_sec} for detailed explanations.

When strong convexity is assumed, Goldstein et al. \cite{Goldstein-2014-fastADM} proposed an $O(1/K^2)$ convergent ADMM for its dual problem. When even more assumptions are made, e.g. the objective function is strongly convex and has Lipschitz continuous gradient, or subdifferentials of the
underlying functions are piecewise linear multifunctions, linear convergence can be obtained \cite{Deng-2012-linearADM,Luo-2012-LinearADM,Giselsson-2016,han-2016,boley-2013}. Some researchers studied the first-order primal-dual algorithm for the saddle-point problem, which includes problem (\ref{problem}) as a special case. For example, Chambolle and Pock \cite{Chambolle-2011} established the $O(1/K)$ ergodic convergence rate for the general convex problems, the accelerated $O(1/K^2)$ convergence rate when the primal or the dual objective is uniformly convex and the linear convergence rate when both are uniformly convex. Chen et al. \cite{Chen-2014} combined Nesterov's second scheme with the primal-dual algorithm and also established the $O(1/K)$ ergodic convergence rate.

\subsection{Contributions}

Although the $O(1/K)$ convergence rate of ADMM and its accelerated versions is widely studied in the literatures, they all need an ergodic averaging \cite{He-2011-Rate,Chen-2015-IPADM,ouyang-2015-fastLADM,Chambolle-2011,Chen-2014}, which may destroy the sparseness and low-rankness in sparse and low-rank learning. As far as we know, there is no literature establishing the $O(1/K)$ nonergodic convergence rate of ADMM type methods for the general convex problem (\ref{problem}). Moreover, as proved in \cite{Davis-2014-DR}, the nonergodic convergence rate of the traditional ADMM is $o(1/\sqrt{K})$ and it will be shown in Section \ref{ADMM-tight} that this rate is tight. In this paper, we aim to give the first $O(1/K)$ nonergodic convergent ADMM type method.

We modify the accelerated ADMM proposed in \cite{ouyang-2015-fastLADM} and give an $O(1/K)$ nonergodic analysis satisfying $|F(\x^K)-F(\x^*)|\leq O(1/K)$ and $\|\A\x^K-\b\|\leq O(1/K)$. Compared with the $O(1/K)$ ergodic rate in \cite{ouyang-2015-fastLADM,He-2011-Rate}, our result is in a nonergodic sense and thus enjoys a more favorable sparseness and low-rankness in applications of sparse and low-rank learning. Compared with the nonergodic rate in \cite{Davis-2014-DR,ouyang-2015-fastLADM}, we improve it from $o(1/\sqrt{K})$  to $O(1/K)$.

We also show that the lower complexity bound of ADMM type methods for the separable linearly constrained convex problems is $O(1/K)$ when each $F_i$ is nonsmooth and non-strongly convex, which means that the convergence rate of ADMM type methods cannot be better than $O(1/K)$ no matter how it is accelerated. Thus our method is optimal.

\section{Review of the Accelerated ADMM in \cite{ouyang-2015-fastLADM}}\label{review_sec}

In this section, we first review the accelerated ADMM in \cite{ouyang-2015-fastLADM} for problem (\ref{problem}), which consists of the following steps:\footnote{We simplify some parameter settings and extend the class of problems it is solving, but the algorithm framework remains the same as \cite{ouyang-2015-fastLADM}.}
\begin{subequations} \begin{align}
&\y_i^k=(1-\theta^k)\x_i^k+\theta^k\z_i^k,i=1,2,\label{step1_LADMM_n2}\\
&\z_1^{k+1}=\argmin_{\z_1} f_1(\z_1^k)+\<\nabla f_1(\y_1^k),\z_1-\z_1^k\>+\frac{\theta^k L_1}{2}\|\z_1-\z_1^k\|^2+h_1(\z_1)\notag\\
&\hspace*{0.1cm}+\hspace*{-0.05cm}\<\hat\lambda^k,\A_1\z_1\>\hspace*{-0.05cm}+\hspace*{-0.05cm}\<\beta\A_1^T(\A_1\z_1^k\hspace*{-0.05cm}+\hspace*{-0.05cm}\A_2\z_2^k\hspace*{-0.05cm}-\hspace*{-0.05cm}\b),\z_1\hspace*{-0.05cm}-\hspace*{-0.05cm}\z_1^k\>\hspace*{-0.05cm}+\hspace*{-0.05cm}\frac{\beta\|\A_1\|^2}{2}\|\z_1\hspace*{-0.05cm}-\hspace*{-0.05cm}\z_1^k\|^2,\label{step2_LADMM_n2}\\
&\z_2^{k+1}=\argmin_{\z_2} f_2(\z_2^k)+\<\nabla f_2(\y_2^k),\z_2-\z_2^k\>+\frac{\theta^k L_2}{2}\|\z_2-\z_2^k\|^2+h_2(\z_2)\notag\\
&\hspace*{0.1cm}+\hspace*{-0.07cm}\<\hat\lambda^k,\A_2\z_2\>\hspace*{-0.07cm}+\hspace*{-0.07cm}\<\beta\A_2^T(\A_1\z_1^{k+1}\hspace*{-0.07cm}+\hspace*{-0.07cm}\A_2\z_2^k\hspace*{-0.07cm}-\hspace*{-0.07cm}\b,\z_2\hspace*{-0.07cm}-\hspace*{-0.07cm}\z_2^k\>\hspace*{-0.07cm}+\hspace*{-0.07cm}\frac{\beta\|\A_2\|^2}{2}\|\z_2\hspace*{-0.07cm}-\hspace*{-0.07cm}\z_2^k\|^2,\label{step3_LADMM_n2}\\
&\x_1^{k+1}=(1-\theta^k)\x_1^k+\theta^k\z_1^{k+1},\label{step41_LADMM_n2}\\
&\x_2^{k+1}=(1-\theta^k)\x_2^k+\theta^k\z_2^{k+1},\label{step42_LADMM_n2}\\
&\hat\lambda^{k+1}=\hat\lambda^k+\beta(\A_1\z_1^{k+1}+\A_2\z_2^{k+1}-\b),\label{step5_LADMM_n2}
\end{align} \end{subequations}
where $\theta^k$ satisfies $\frac{1}{(\theta^{k-1})^2}\geq\frac{1-\theta^k}{(\theta^k)^2}$. Since the regularizer $h(\x)$ acts directly on $\z$ in (\ref{step2_LADMM_n2})-(\ref{step3_LADMM_n2}), $(\z_1^{k+1},\z_2^{k+1})$ has the property promoted by $h(\x)$ and the convergence measured at $(\z_1^{K},\z_2^K)$ is in the nonergodic sense. In fact, in sparse or low-rank learning, we often use the $l_1$- norm and nuclear norm as the regularization. The proximal operator of the $l_1$-norm is the soft-thresholding \cite{denoho-tit}, which is defined as
\begin{eqnarray}
\begin{aligned}\notag
\argmin_{\z}\|\z\|_1+\frac{\gamma}{2}\|\z-\w\|^2=\left\{
      \begin{array}{lcl}
        \w_i-\frac{1}{\gamma}, && \mbox{if } \w_i\geq \frac{1}{\gamma},  \\
        \w_i+\frac{1}{\gamma}, && \mbox{if } \w_i\leq -\frac{1}{\gamma},  \\
        0, && \mbox{otherwise.}
      \end{array}
    \right.
\end{aligned}
\end{eqnarray}
Thus, if we use $h(\z)=\|\z\|_1$, $(\z_1^{k+1},\z_2^{k+1})$ tends to be sparse during the iterations. Similarly, the proximal operator of the nuclear norm is the singular value thresholding \cite{cai-2010-siam}, which is defined as
\begin{eqnarray}
\begin{aligned}\notag
\argmin_{\Z}\|\Z\|_*+\frac{\gamma}{2}\|\Z-\W\|_F^2=\U\hat\Sigma\V^T,
\end{aligned}
\end{eqnarray}
where we let $\U\Sigma\V^T=\W$ be its SVD and
    \begin{eqnarray}
    \begin{aligned}\notag
    \hat\Sigma_{i,i}=\left\{
      \begin{array}{lcl}
        \Sigma_{i,i}-\frac{1}{\gamma}, && \mbox{if } \Sigma_{i,i}\geq \frac{1}{\gamma},  \\
        \Sigma_{i,i}+\frac{1}{\gamma}, && \mbox{if } \Sigma_{i,i}\leq -\frac{1}{\gamma},  \\
        0, && \mbox{otherwise.}
      \end{array}
    \right.
    \end{aligned}
    \end{eqnarray}
Thus, if we use $h(\Z)=\|\Z\|_*$, $(\Z_1^{k+1},\Z_2^{k+1})$ tends to be low-rank during the iterations.

Accordingly, $\x^K$ is a convex combination of $\z^1,\cdots,\z^K$: $\x^K=\frac{1}{\sum_{k=1}^K \frac{1}{\theta^{k-1}}}\sum_{k=1}^K \frac{\z^k}{\theta^{k-1}}$ and so it is an ergodic result measured at $(\x_1^{K},\x_2^K)$. The zeros may lie in different positions of $\z^1,\cdots,\z^K$ for sparse learning (or in different positions of their singular values for low-rank learning) and thus $\x^K$ may not be sparse or low-rank any more. It is proved in \cite{ouyang-2015-fastLADM} that (\ref{step1_LADMM_n2})-(\ref{step5_LADMM_n2}) has the $O(1/K)$ ergodic convergence rate measured at $(\x_1^{K},\x_2^K)$.

We can see that the accelerated ADMM in \cite{ouyang-2015-fastLADM} is a direct combination of Nesterov's second acceleration scheme and the traditional LADMM. Nesterov's acceleration scheme only influences the linearization of $f_i$ and cannot improve the convergence rate of the traditional ADMM. In fact, we can consider the special case of $f_i(\x_i)=0,i=1,2$ (correspondingly, $L_i=0$) and omit the linearization of the augmented term (or let $\A_1=\A_2=\I$ for simplicity). In this case, procedure (\ref{step1_LADMM_n2})-(\ref{step5_LADMM_n2}) reduces to:
\begin{subequations} \begin{align}
&\z_1^{k+1}=\argmin_{\z_1} h_1(\z_1)+\<\hat\lambda^k,\A_1\z_1\>+\frac{\beta}{2}\|\A_1\z_1+\A_2\z_2^k-\b\|^2,\label{step1_ADMM_n2}\\
&\z_2^{k+1}=\argmin_{\z_2} h_2(\z_2)+\<\hat\lambda^k,\A_2\z_2\>+\frac{\beta}{2}\|\A_1\z_1^{k+1}+\A_2\z_2-\b\|^2,\label{step2_ADMM_n2}\\
&\x_1^{k+1}=(1-\theta^k)\x_1^k+\theta^k\z_1^{k+1},\label{step3_ADMM_n2}\\
&\x_2^{k+1}=(1-\theta^k)\x_2^k+\theta^k\z_2^{k+1},\label{step4_ADMM_n2}\\
&\hat\lambda^{k+1}=\hat\lambda^k+\beta(\A_1\z_1^{k+1}+\A_2\z_2^{k+1}-\b).\label{step5_ADMM_n2}
\end{align} \end{subequations}
We can see that procedure (\ref{step1_ADMM_n2})-(\ref{step5_ADMM_n2}) reduces to the traditional ADMM and (\ref{step3_ADMM_n2})-(\ref{step4_ADMM_n2}) has no influence on the iterations of the traditional ADMM. It only gives a different way of ergodic averaging. Thus the nonergodic convergence rate of procedure (\ref{step1_ADMM_n2})-(\ref{step5_ADMM_n2}) measured at $(\z_1^{K},\z_2^K)$ remains $o(1/\sqrt{K})$. Since (\ref{step1_ADMM_n2})-(\ref{step5_ADMM_n2}) is a special case of (\ref{step1_LADMM_n2})-(\ref{step5_LADMM_n2}), we can have that the nonergodic rate of procedure (\ref{step1_LADMM_n2})-(\ref{step5_LADMM_n2}) measured at $(\z_1^{K},\z_2^K)$ should not be better than $o(1/\sqrt{K})$.
\section{ALADMM-NE with $O(1/K)$ Nonergodic Convergence Rate}\label{general_convex_section}
In this section, we give our Accelerated LADMM with NonErgodic convergence rate (ALADMM-NE). We first provide an equivalent description of (\ref{step1_LADMM_n2})-(\ref{step5_LADMM_n2}) for the smooth case of problem (\ref{problem}) in Section \ref{equi_sec}, which motivates our nonergodic algorithm for the nonsmooth case in Section \ref{alg_sec}. Then we give the convergence rate analysis in Section \ref{rate_sec} and at last, we discuss the advantage and disadvantage of the accelerated ADMM in Section \ref{tip_sec}.
\subsection{An Equivalent Algorithm for the Smooth Problem}\label{equi_sec}
In this section, we give an equivalent description of (\ref{step1_LADMM_n2})-(\ref{step5_LADMM_n2}) for the smooth case of problem (\ref{problem}) with $h_i(\x)=0,i=1,2$:
\begin{subequations} \begin{align}
&\y_i^k=\x_i^k+\frac{\theta^k(1-\theta^{k-1})}{\theta^{k-1}}(\x_i^k-\x_i^{k-1}),i=1,2\label{step1_LADMM_sn1},\\
&\x_1^{k+1}=\argmin_{\x_1} f_1(\y_1^k)+\<\nabla f_1(\y_1^k),\x_1-\y_1\>+\frac{L_1}{2}\|\x_1-\y_1^k\|^2+\<\lambda^k,\A_1\x_1\>\notag\\
&\qquad\quad+\frac{\beta}{\theta^k}\<\A_1^T(\A_1\y_1^k+\A_2\y_2^k-\b),\x_1-\y_1^k\>+\frac{\beta\|\A_1\|^2}{2\theta^k}\|\x_1-\y_1^k\|^2\label{step2_LADMM_sn1},\\
&\x_2^{k+1}=\argmin_{\x_2} f_2(\y_2^k)+\<\nabla f_2(\y_2^k),\x_2-\y_2\>+\frac{L_2}{2}\|\x_2-\y_2^k\|^2+\<\lambda^k,\A_2\x_2\>\notag\\
&\quad\quad+\frac{\beta}{\theta^k}\<\A_2^T(\A_1\x_1^{k+1}+\A_2\y_2^k-\b),\x_2-\y_2^k\>+\frac{\beta\|\A_2\|^2}{2\theta^k}\|\x_2-\y_2^k\|^2\label{step3_LADMM_sn1},\\
&\lambda^{k+1}=\lambda^k+\beta\tau(\A\x^{k+1}-\b)\label{step4_LADMM_sn1},
\end{align} \end{subequations}
for some $1>\tau>0.5$, $\theta^0=1$ and $\theta^{k+1}=\frac{1}{1-\tau+\frac{1}{\theta^k}}$, which leads to $\frac{1}{(\theta^{k-1})^2}\geq\frac{1-\theta^k}{(\theta^k)^2}$ and thus coincides with the requirement for (\ref{step2_LADMM_n2})-(\ref{step5_LADMM_n2}). It can be observed that if we set $\tau=1$, then $\theta^k=1,\y_i^k=\x_i^k,\forall k$, and (\ref{step1_LADMM_sn1})-(\ref{step4_LADMM_sn1}) reduces to the traditional LADMM. At first glance, (\ref{step1_LADMM_n2})-(\ref{step5_LADMM_n2}) combines ADMM with Nesterov's second acceleration scheme while (\ref{step1_LADMM_sn1})-(\ref{step4_LADMM_sn1}) uses Nesterov's first acceleration scheme.
\begin{proposition}
The sequence $(\x_1^k,\x_2^k)$ produced in (\ref{step1_LADMM_n2})-(\ref{step5_LADMM_n2}) and (\ref{step1_LADMM_sn1})-(\ref{step4_LADMM_sn1}) are equivalent when $h_i(\x)=0,i=1,2$.
\end{proposition}
\begin{proof}
We derive each step of (\ref{step1_LADMM_sn1})-(\ref{step4_LADMM_sn1}) from (\ref{step1_LADMM_n2})-(\ref{step5_LADMM_n2}). From (\ref{step1_LADMM_n2}), (\ref{step41_LADMM_n2}) and (\ref{step42_LADMM_n2}), we have
\begin{subequations} \begin{align}
\y_i^k&=(1-\theta^k)\x_i^k+\theta^k\z_i^k\notag\\
&=(1-\theta^k)\x_i^k+\frac{\theta^k}{\theta^{k-1}}\left(\x_i^{k}-(1-\theta^{k-1})\x_i^{k-1}\right)\notag\\
&=\x_i^k+\frac{\theta^k(1-\theta^{k-1})}{\theta^{k-1}}(\x_i^k-\x_i^{k-1}),\notag
\end{align} \end{subequations}
which is (\ref{step1_LADMM_sn1}). From the optimality condition of (\ref{step2_LADMM_n2}), we have
\begin{subequations} \begin{align}
0&=\nabla f_1(\y_1^k)+\theta^k L_1(\z_1^{k+1}-\z_1^k)+\A_1^T\hat\lambda^k+\beta\A_1^T(\A_1\z_1^k+\A_2\z_2^k-\b)\notag\\
&\qquad\quad+\beta\|\A_1\|^2(\z_1^{k+1}-\z_1^k)\notag\\
&\overset{\ref{step1_LADMM_n2},\ref{step41_LADMM_n2}}{=}\nabla f_1(\y_1^k)+ L_1(\x_1^{k+1}-\y_1^k)+\A_1^T\hat\lambda^k+\beta\A_1^T(\A_1\z_1^k+\A_2\z_2^k-\b)\notag\\
&\qquad\quad+\frac{\beta\|\A_1\|^2}{\theta^k}(\x_1^{k+1}-\y_1^k)\notag\\
&\overset{\ref{step1_LADMM_n2}}{=}\nabla f_1(\y_1^k)+ L_1(\x_1^{k+1}-\y_1^k)+\A_1^T\hat\lambda^k\notag\\
&\qquad\quad+\beta\A_1^T\left(\frac{\A_1\y_1^k}{\theta^k}-\frac{1-\theta^k}{\theta^k}\A_1\x_1^k+\frac{\A_2\y_2^k}{\theta^k}-\frac{1-\theta^k}{\theta^k}\A_2\x_2^k-\b\right)\notag\\
&\qquad\quad+\frac{\beta\|\A_1\|^2}{\theta^k}(\x_1^{k+1}-\y_1^k)\notag\\
&=\nabla f_1(\y_1^k)+ L_1(\x_1^{k+1}-\y_1^k)+\A_1^T\hat\lambda^k-\frac{\beta(1-\theta^k)}{\theta^k}\A_1^T(\A_1\x_1^k+\A_2\x_2^k-\b)\notag\\
&\qquad\quad+\frac{\beta}{\theta^k}\A_1^T\left(\A_1\y_1^k+\A_2\y_2^k-\b\right)+\frac{\beta\|\A_1\|^2}{\theta^k}(\x_1^{k+1}-\y_1^k)\notag\\
&=\nabla f_1(\y_1^k)+ L_1(\x_1^{k+1}-\y_1^k)+\A_1^T\lambda^k+\frac{\beta}{\theta^k}\A_1^T\left(\A_1\y_1^k+\A_2\y_2^k-\b\right)\notag\\
&\qquad\quad+\frac{\beta\|\A_1\|^2}{\theta^k}(\x_1^{k+1}-\y_1^k),\notag
\end{align} \end{subequations}
where we define $\lambda^k=\hat\lambda^k-\frac{\beta(1-\theta^k)}{\theta^k}(\A_1\x_1^k+\A_2\x_2^k-\b)$. It is exactly the optimality condition of (\ref{step2_LADMM_sn1}). Similarly, from the optimality condition of (\ref{step3_LADMM_n2}), we also have
\begin{subequations} \begin{align}
0&=\nabla f_2(\y_2^k)+\theta^k L_2(\z_2^{k+1}-\z_2^k)+\A_2^T\hat\lambda^k+\beta\A_2^T(\A_1\z_1^{k+1}+\A_2\z_2^k-\b)\notag\\
&\qquad\quad+\beta\|\A_2\|^2(\z_2^{k+1}-\z_2^k)\notag\\
&=\nabla f_2(\y_2^k)+ L_2(\x_2^{k+1}-\y_2^k)+\A_2^T\hat\lambda^k+\beta\A_2^T(\A_1\z_1^{k+1}+\A_2\z_2^k-\b)\notag\\
&\qquad\quad+\frac{\beta\|\A_2\|^2}{\theta^k}(\x_2^{k+1}-\y_2^k)\notag\\
&\overset{\ref{step41_LADMM_n2},\ref{step1_LADMM_n2}}{=}\nabla f_2(\y_2^k)+ L_2(\x_2^{k+1}-\y_2^k)+\A_2^T\hat\lambda^k\notag\\
&\qquad\quad+\beta\A_2^T\left(\frac{\A_1\x_1^{k+1}}{\theta^k}-\frac{1-\theta^k}{\theta^k}\A_1\x_1^k+\frac{\A_2\y_2^k}{\theta^k}-\frac{1-\theta^k}{\theta^k}\A_2\x_2^k-\b\right)\notag\\
&\qquad\quad+\frac{\beta\|\A_2\|^2}{\theta^k}(\x_2^{k+1}-\y_2^k)\notag
\end{align} \end{subequations}
\begin{subequations} \begin{align}
&=\nabla f_2(\y_2^k)+ L_2(\x_2^{k+1}-\y_2^k)+\A_2^T\lambda^k+\frac{\beta}{\theta^k}\A_2^T\left(\A_1\x_1^{k+1}+\A_2\y_2^k-\b\right)\notag\\
&\qquad\quad+\frac{\beta\|\A_2\|^2}{\theta^k}(\x_2^{k+1}-\y_2^k),\notag
\end{align} \end{subequations}
which is the optimality condition of (\ref{step3_LADMM_sn1}). From the definition of $\lambda^k$, we have
\begin{subequations} \begin{align}
&\lambda^{k+1}-\lambda^k\notag\\
&=\hat\lambda^{k+1}-\hat\lambda^k-\frac{\beta(1-\theta^{k+1})}{\theta^{k+1}}(\A\x^{k+1}-\b)+\frac{\beta(1-\theta^k)}{\theta^k}(\A\x^k-\b)\notag\\
&\overset{\ref{step5_LADMM_n2}}{=}\beta(\A\z^{k+1}-\b)-\frac{\beta(1-\theta^{k+1})}{\theta^{k+1}}(\A\x^{k+1}-\b)+\frac{\beta(1-\theta^k)}{\theta^k}(\A\x^k-\b)\notag\\
&\overset{\ref{step41_LADMM_n2},\ref{step42_LADMM_n2}}{=}\beta\left(\frac{\A\x^{k+1}-(1-\theta^k)\A\x^k}{\theta^k}-\b\right)-\frac{\beta(1-\theta^{k+1})}{\theta^{k+1}}(\A\x^{k+1}-\b)\notag\\
&\qquad\quad+\frac{\beta(1-\theta^k)}{\theta^k}(\A\x^k-\b)\notag\\
&=\beta\left(\frac{(\A\x^{k+1}-\b)-(1-\theta^k)(\A\x^k-\b)}{\theta^k}\right)-\frac{\beta(1-\theta^{k+1})}{\theta^{k+1}}(\A\x^{k+1}-\b)\notag\\
&\qquad\quad+\frac{\beta(1-\theta^k)}{\theta^k}(\A\x^k-\b)\notag\\
&=\beta\left(\frac{1}{\theta^k}-\frac{1-\theta^{k+1}}{\theta^{k+1}}\right)(\A\x^{k+1}-\b)\notag\\
&=\beta\tau(\A\x^{k+1}-\b),\notag
\end{align} \end{subequations}
where we define $\tau=\frac{1}{\theta^k}-\frac{1-\theta^{k+1}}{\theta^{k+1}}$ and it is the same with (\ref{step4_LADMM_sn1}).$\hfill\blacksquare$
\end{proof}
\subsection{The Nonergodic Algorithm for the Nonsmooth Problem}\label{alg_sec}

\begin{algorithm}[tb]
   \caption{Accelerated LADMM with NonErgodic convergence rate (ALADMM-NE)}
   \label{gc_serial_alg}
\begin{algorithmic}
   \STATE Initialize $\lambda^0$, $\x_i^0=\x_i^{-1},i=1,2$, $1>\tau> 0.5$, $\beta>0$, $\theta^0=1$, $\theta^{-1}=1/\tau$.
   \FOR{$k=0,1,2,\cdots$}
   \STATE Update $\y_i^k,i=1,2$ using (\ref{step1_LADMM_sn1}),\\
   \STATE Update $\x_1^{k+1}$ and $\x_2^{k+1}$ serially, using (\ref{step2_LADMM_n1}) and (\ref{step3_LADMM_n1}), respectively,\\
   \STATE Update $\lambda^{k+1}$ using (\ref{step4_LADMM_sn1}),\\
   $\theta^{k+1}=\frac{1}{1-\tau+\frac{1}{\theta^k}}$.
   \ENDFOR
\end{algorithmic}
\end{algorithm}

From the discussion in Section \ref{review_sec}, we know that the accelerated ADMM proposed in \cite{ouyang-2015-fastLADM} has the $o(1/\sqrt{K})$ nonergodic convergence rate measured at $(\z_1^K,\z_2^k)$ and the $O(1/K)$ ergodic convergence rate measured at $(\x_1^k,\x_2^k)$. We want to have an algorithm with the faster $O(1/K)$ nonergodic convergence rate. After establishing the equivalence between (\ref{step1_LADMM_n2})-(\ref{step5_LADMM_n2}) and (\ref{step1_LADMM_sn1})-(\ref{step4_LADMM_sn1}), a straightforward intuition is to put the nonsmooth term $h_i(\x)$ in steps (\ref{step2_LADMM_sn1}) and (\ref{step3_LADMM_sn1}) directly:
\begin{subequations} \begin{align}
\x_1^{k+1}=&\argmin_{\x_1} f_1(\y_1^k)+\<\nabla f_1(\y_1^k),\x_1-\y_1\>+\frac{L_1}{2}\|\x_1-\y_1^k\|^2+h_1(\x_1)\notag\\
&+\hspace*{-0.1cm}\<\lambda^k,\A_1\x_1\>\hspace*{-0.1cm}+\hspace*{-0.1cm}\frac{\beta}{\theta^k}\<\A_1^T(\A_1\y_1^k\hspace*{-0.1cm}+\hspace*{-0.1cm}\A_2\y_2^k\hspace*{-0.1cm}-\hspace*{-0.1cm}\b),\x_1\hspace*{-0.1cm}-\hspace*{-0.1cm}\y_1^k\>\hspace*{-0.1cm}+\hspace*{-0.1cm}\frac{\beta\|\A_1\|^2}{2\theta^k}\|\x_1\hspace*{-0.1cm}-\hspace*{-0.1cm}\y_1^k\|^2\label{step2_LADMM_n1},\\
\x_2^{k+1}=&\argmin_{\x_2} f_2(\y_2^k)+\<\nabla f_2(\y_2^k),\x_2-\y_2\>+\frac{L_2}{2}\|\x_2-\y_2^k\|^2+h_2(\x_2)\notag\\
&+\hspace*{-0.1cm}\<\hspace*{-0.05cm}\lambda^k\hspace*{-0.05cm},\hspace*{-0.05cm}\A_2\x_2\>\hspace*{-0.1cm}+\hspace*{-0.1cm}\frac{\beta}{\theta^k}\hspace*{-0.1cm}\<\hspace*{-0.05cm}\A_2^T(\A_1\x_1^{k+1}\hspace*{-0.1cm}+\hspace*{-0.1cm}\A_2\y_2^k\hspace*{-0.1cm}-\hspace*{-0.1cm}\b),\x_2\hspace*{-0.1cm}-\hspace*{-0.1cm}\y_2^k\>\hspace*{-0.1cm}+\hspace*{-0.1cm}\frac{\beta\|\A_2\|^2}{2\theta^k}\|\x_2\hspace*{-0.1cm}-\hspace*{-0.1cm}\y_2^k\|^2\label{step3_LADMM_n1}.
\end{align} \end{subequations}

We describe the new method in Algorithm \ref{gc_serial_alg}. Due to the different positions of the terms $h_i(\x_i),i=1,2$, Algorithm \ref{gc_serial_alg} and procedure (\ref{step1_LADMM_n2})-(\ref{step5_LADMM_n2}) are no longer equivalent for the nonsmooth problem. In fact, when we consider the simple case of $f_1(\x_1)=0$ and $f_2(\x_2)=0$ and omit the linearization of the augmented term, (\ref{step1_LADMM_n2})-(\ref{step5_LADMM_n2}) reduces to the traditional ADMM, while Algorithm \ref{gc_serial_alg} reduces to the following iterates:
\begin{subequations} \begin{align}
    &\y_i^k=\x_i^k+\frac{\theta^k(1-\theta^{k-1})}{\theta^{k-1}}(\x_i^k-\x_i^{k-1}),i=1,2\label{acc1},\\
    &\x_1^{k+1}=\argmin_{\x_1}h_1(\x_1)+\<\lambda^k,\A_1\x_1\>+\frac{\beta}{2\theta^k}\|\A_1\x_1+\A_2\y_2^k-\b\|^2,\label{acc2}\\
    &\x_2^{k+1}=\argmin_{\x_2}h_2(\x_2)+\<\lambda^k,\A_2\x_2\>+\frac{\beta}{2\theta^k}\|\A_1\x_1^{k+1}+\A_2\x_2-\b\|^2,\label{acc3}\\
    &\lambda^{k+1}=\lambda^k+\beta\tau(\A\x^{k+1}-\b).\label{acc4}
\end{align} \end{subequations}
We can see that procedure (\ref{acc1})-(\ref{acc4}) is totally different from the traditional ADMM, which verifies that Algorithm \ref{gc_serial_alg} is different from procedure (\ref{step1_LADMM_n2})-(\ref{step5_LADMM_n2}). The analysis in this paper can be easily used to establish the $O(1/K)$ nonergodic convergence rate of procedure (\ref{acc1})-(\ref{acc4}) measured at $\{\x_1^k,\x_2^k\}$. We only consider the complex case of Algorithm \ref{gc_serial_alg} and omit the proof for the simple case of (\ref{acc1})-(\ref{acc4}).

In Algorithm \ref{gc_serial_alg}, $h_i(\x_i)$ acts on $\x_i$ directly and thus it has the property promoted by $h_i(\x_i)$, such as the sparseness or low-rankness if $h_i(\x_i)$ is a sparse or low rank regularizer. So the convergence rate measured at $\x^K$ in Algorithm \ref{gc_serial_alg} is in the nonergodic sense. As comparison, (\ref{step1_LADMM_n2})-(\ref{step5_LADMM_n2}) promotes the sparseness and low-rankness on $\z_i$, and $\x^K$ is a convex combination of $\z^1,\cdots,\z^K$ and it may not be sparse or low-rank any more due to the ergodic averaging. In applications where sparseness or low-rankness is strongly required, we should use the nonergodic solutions and Algorithm \ref{gc_serial_alg} is superior to (\ref{step1_LADMM_n2})-(\ref{step5_LADMM_n2}), since the nonergodic solution in Algorithm \ref{gc_serial_alg} has a faster convergence rate than the nonergodic solution in procedure (\ref{step1_LADMM_n2})-(\ref{step5_LADMM_n2}). We demonstrate the differences in Table \ref{com_table1}. It should be noted that for the smooth case, since $h(\x)$ vanishes, we do not distinguish the ergodic and the nonergodic rates between (\ref{step1_LADMM_n2})-(\ref{step5_LADMM_n2}) and (\ref{step1_LADMM_sn1})-(\ref{step4_LADMM_sn1}).

\begin{table}
\begin{center}
\begin{tabular}{l|c|c|c|c}
\hline\hline
Algorithm & Solution & Ergodic or Nonergodic & Sparse/Low-rank & Convergence Rate \\
\hline\hline
\multirow{2}{*}{(\ref{adm1})-(\ref{adm3})} & $\frac{\sum_{k=1}^K\x^k}{K}$ & Ergodic & No & $O\left(1/K\right)$\\\cline{2-5}
                                               & $\x^K$ & Nonergodic & Yes & $O\left(1/\sqrt{K}\right)$ \\
\hline\hline
\multirow{2}{*}{(\ref{step1_LADMM_n2})-(\ref{step5_LADMM_n2})} & $\x^K$ & Ergodic & No & $O\left(1/K\right)$\\\cline{2-5}
                                               & $\z^K$ & Nonergodic & Yes & $O\left(1/\sqrt{K}\right)$ \\
\hline\hline
Algorithm \ref{gc_serial_alg} & $\x^K$ & Nonergodic & Yes & $O\left(1/K\right)$\\
\hline\hline
\end{tabular}
\end{center}
\caption{Comparing Algorithm \ref{gc_serial_alg} with the original ADMM, (\ref{adm1})-(\ref{adm3}) and the accelerated ADMM, (\ref{step1_LADMM_n2})-(\ref{step5_LADMM_n2}) on the properties of the ergodic and nonergodic solutions.}
\label{com_table1}
\end{table}

\subsection{The Convergence Rate Analysis}\label{rate_sec}

In this section, we prove the $O(1/K)$ convergence rate measured at $\x^K$ for Algorithm \ref{gc_serial_alg}. Due to the different positions of the nonsmooth terms $h_i(\x_i)$, the proof technique for procedure (\ref{step1_LADMM_n2})-(\ref{step5_LADMM_n2}) in \cite{ouyang-2015-fastLADM} cannot be extended to Algorithm \ref{gc_serial_alg} and more efforts are needed for the analysis on Algorithm \ref{gc_serial_alg}. Moreover, Ouyang et al. \cite{ouyang-2015-fastLADM} need the assumption that the primal and dual variables are bounded in order to accomplish the proof. As comparison, we do not need this assumption. This verifies that our proof is totally different from \cite{ouyang-2015-fastLADM}.

ALADMM-NE is an extension of Nesterov's first acceleration scheme from unconstrained problems to constrained ones. For unconstrained problems, a crucial property of Nesterov's first acceleration scheme is
\begin{equation}
\frac{F(\x^{k+1})-F(\x^*)}{(\theta^k)^2}-\frac{F(\x^k)-F(\x^*)}{(\theta^{k-1})^2}\leq \delta\left(\|\z^k-\x^*\|^2-\|\z^{k+1}-\x^*\|^2\right).\label{APG-descent}
\end{equation}
The main step in the convergence rate proof of ALADMM-NE is to construct a counterpart of (\ref{APG-descent}) for both the objective and the constraint functions. Proposition \ref{serial_gc_general_convex_decent_lemma} plays such a role for the objective. As comparison, the traditional ADMM \cite{Lin-2015-LADMPSAP} can prove a similar result in the form of
\begin{eqnarray}
&&F(\x^k)-F(\x^*)+\<\lambda^*,\A\x^k-\b\>\notag\\
&\leq& \delta\left(\|\x^k-\x^*\|^2-\|\x^{k+1}-\x^*\|^2\right)+\kappa\left(\|\lambda^k-\lambda^*\|^2-\|\lambda^{k+1}-\lambda^*\|^2\right),\notag
\end{eqnarray}
which can only lead to the ergodic result after telescoping.

\begin{proposition}\label{serial_gc_general_convex_decent_lemma}
Assume that $f_i(\x_i)$ is convex and $L_i$-smooth, $h_i(\x_i)$ is convex, i=1,2. Let $\frac{1-\theta^k}{\theta^k}=\frac{1}{\theta^{k-1}}-\tau$ with $0< \tau< 1$ and $\theta^0=1$. For Algorithm \ref{gc_serial_alg}, we have
\begin{subequations} \begin{align}
&\frac{1}{\theta^k}\left(F(\x^{k+1})\hspace*{-0.1cm}-\hspace*{-0.1cm}F(\x^*)\hspace*{-0.1cm}+\hspace*{-0.1cm}\<\lambda^*,\A\x^{k+1}\hspace*{-0.1cm}-\hspace*{-0.1cm}\b\>\right)\hspace*{-0.1cm}-\hspace*{-0.1cm}\frac{1}{\theta^{k-1}}\left( F(\x^k)\hspace*{-0.1cm}-\hspace*{-0.1cm}F(\x^*)\hspace*{-0.1cm}+\hspace*{-0.1cm}\<\lambda^*,\A\x^k\hspace*{-0.1cm}-\hspace*{-0.1cm}\b\> \right)\notag\\
&+\tau\left( F(\x^k)-F(\x^*)+\<\lambda^*,\A\x^k-\b\> \right)\notag\\
\leq&\frac{1}{2\beta}\left(\|\hat\lambda^k-\lambda^*\|^2-\|\hat\lambda^{k+1}-\lambda^*\|^2\right)+\frac{\eta_2^k}{2}\|\d_2^k-\x_2^*\|^2-\frac{\eta_2^{k+1}}{2}\|\d_2^{k+1}-\x_2^*\|^2\notag\\
&+\left(\frac{\eta_1^k}{2}\|\d_1^k-\x_1^*\|^2-\frac{\beta}{2}\|\A_1\d_1^k-\A_1\x^*\|^2\right)\notag\\
&-\left(\frac{\eta_1^{k+1}}{2}\|\d_1^{k+1}-\x_1^*\|^2-\frac{\beta}{2}\|\A_1\d_1^{k+1}-\A_1\x^*\|^2\right),\notag
\end{align} \end{subequations}
where $\hat\lambda^k=\lambda^k+\frac{\beta(1-\theta^k)}{\theta^k}\left(\A\x^k-\b\right)$, $\eta_i^k=L_i\theta^k+\beta\|\A_i\|_2^2 $, $\d_i^{k+1}=\frac{\x_i^{k+1}}{\theta^k}-\frac{1-\theta^k}{\theta^k}\x_i^k$, $\d_i^0=\x_i^0,i=1,2,$ and $\{\x^*,\lambda^*\}$ is any KKT point.
\end{proposition}

Before proving Proposition \ref{serial_gc_general_convex_decent_lemma}, we first prove the following Lemma.
\begin{lemma}\label{lemma_lambda}
Let $\overline\lambda_2^{k+1}=\lambda^k+\frac{\beta}{\theta^k}\left( \A_1\x_1^{k+1}+\A_2\y_2^k-\b\right)$. Then for Algorithm \ref{gc_serial_alg}, we have
\begin{subequations} \begin{align}
\begin{aligned}
&\theta^k\b+(1-\theta^k)\A\x^k-\A\x^{k+1}=\frac{\theta^k}{\beta}\left(\hat\lambda^k-\hat\lambda^{k+1}\right),\notag\\
&\frac{\theta^k}{2\beta}\|\hat\lambda^{k+1}-\overline\lambda_2^{k+1}\|^2\leq \frac{\beta}{2\theta^k}\|\A_2\|_2^2\|\x_2^{k+1}-\y_2^k\|^2,\notag\\
&\hat\lambda^{K+1}-\hat\lambda^0=\sum_{k=0}^K \left[ \beta\frac{\A\x^{k+1}-\b}{\theta^k}-\beta\frac{\A\x^{k}-\b}{\theta^{k-1}}+\beta\tau\left(\A\x^{k}-\b\right)\right].\notag
\end{aligned}
\end{align} \end{subequations}
\end{lemma}
\begin{proof}
From $\hat\lambda^k=\lambda^k+\frac{\beta(1-\theta^k)}{\theta^k}\left(\A\x^k-\b\right)$, $\frac{1-\theta^{k+1}}{\theta^{k+1}}=\frac{1}{\theta^k}-\tau$ and $\lambda^{k+1}=\lambda^k+\beta\tau\left( \sum_{i=1}^2\A_i\x_i^{k+1}-\b\right)$ we have
\begin{subequations} \begin{align}
&\hat\lambda^{k+1}=\lambda^{k+1}+\beta\frac{(1-\theta^{k+1})\left(\sum_{i=1}^2\A_i\x_i^{k+1}-\b\right)}{\theta^{k+1}}\notag\\
=&\lambda^{k+1}+\beta\left(\frac{1}{\theta^k}-\tau\right)\left(\sum_{i=1}^2\A_i\x_i^{k+1}-\b\right)\notag\\
=&\lambda^k+\beta\tau\left(\sum_{i=1}^2\A_i\x_i^{k+1}-\b\right)+\beta\left(\frac{1}{\theta^k}-\tau\right)\left(\sum_{i=1}^2\A_i\x_i^{k+1}-\b\right)\notag\\
=&\lambda^k+\frac{\beta}{\theta^k}\left(\sum_{i=1}^2\A_i\x_i^{k+1}-\b\right)\label{serial_gc_hat_lambda}\\
=&\hat\lambda^k-\beta\frac{(1-\theta^k)\left(\sum_{i=1}^2\A_i\x_i^k-\b\right)}{\theta^{k}}+\frac{\beta}{\theta^k}\left(\sum_{i=1}^2\A_i\x_i^{k+1}-\b\right)\label{serial_gc_dif_hat_lambda1}\\
=&\hat\lambda^k-\frac{\beta}{\theta^k}\left(\theta^k\b+(1-\theta^k)\sum_{i=1}^2\A_i\x_i^k-\sum_{i=1}^2\A_i\x_i^{k+1}\right).\notag
\end{align} \end{subequations}
On the other hand, from (\ref{serial_gc_hat_lambda}) and the definition of $\overline\lambda_2^{k+1}$ we have
\begin{subequations} \begin{align}
\frac{\theta^k}{2\beta}\|\hat\lambda^{k+1}-\overline\lambda_2^{k+1}\|^2=&\frac{\theta^k}{2\beta}\left\|\frac{\beta}{\theta^k}\A_2(\x_2^{k+1}-\y_2^k) \right\|^2\leq \frac{\beta}{2\theta^k}\|\A_2\|_2^2\|\x_2^{k+1}-\y_2^k\|^2.\notag
\end{align} \end{subequations}
From (\ref{serial_gc_dif_hat_lambda1}) and $\frac{1-\theta^k}{\theta^k}=\frac{1}{\theta^{k-1}}-\tau$ we have
\begin{subequations} \begin{align}
&\hspace*{-2cm}\hat\lambda^{K+1}-\hat\lambda^0\notag\\
&\hspace*{-2.5cm}=\sum_{k=0}^K\left(\hat\lambda^{k+1}-\hat\lambda^k\right)\notag\\
&\hspace*{-2.5cm}=\sum_{k=0}^K \left[\beta\frac{\sum_{i=1}^2\A_i\x_i^{k+1}-\b}{\theta^k}-\beta\frac{1-\theta^k}{\theta^k}\left(\sum_{i=1}^2\A_i\x_i^{k}-\b\right)\right]\notag
\end{align} \end{subequations}
\begin{subequations} \begin{align}
=&\sum_{k=0}^K \left[ \beta\frac{\sum_{i=1}^2\A_i\x_i^{k+1}-\b}{\theta^k}-\beta\frac{\sum_{i=1}^2\A_i\x_i^{k}-\b}{\theta^{k-1}}+\beta\tau\left(\sum_{i=1}^2\A_i\x_i^{k}-\b\right)\right]\notag.
\end{align} \end{subequations}
$\hfill\blacksquare$
\end{proof}

Then we can prove Proposition \ref{serial_gc_general_convex_decent_lemma} using Lemma \ref{lemma_lambda}.
\begin{proof}
Let $\overline\lambda_1^{k+1}=\lambda^k+\frac{\beta}{\theta^k}\left( \sum_{i=1}^2\A_i\y_i^k-\b\right)$. From the optimality conditions of (\ref{step2_LADMM_n1}) and (\ref{step3_LADMM_n1}), we have
\begin{subequations} \begin{align}
0\in& \nabla f_i(\y_i^k)+\partial h_i(\x_i^{k+1})+\A_i^T\overline\lambda_i^{k+1}+\left(L_i+\frac{\beta\|\A_i\|_2^2 }{\theta^k}\right)(\x_i^{k+1}-\y_i^k),\notag
\end{align} \end{subequations}
From the convexity of $h_i(\x_i)$ we have
\begin{subequations} \begin{align}
&h_i(\x_i)-h_i(\x_i^{k+1})\notag\\
\geq& -\<\nabla f_i(\y_i^k)+\A_i^T\overline\lambda_i^{k+1}+\left(L_i+\frac{\beta\|\A_i\|_2^2 }{\theta^k}\right)(\x_i^{k+1}-\y_i^k),\x_i-\x_i^{k+1}\>.\notag
\end{align} \end{subequations}
On the other hand, since $f_i$ is $L_i$-smooth and convex, we have
\begin{subequations} \begin{align}
f_i(\x_i^{k+1})\leq&f_i(\y_i^k)+\<\nabla f_i(\y_i^k),\x_i^{k+1}-\y_i^k\>+\frac{L_i}{2}\|\x_i^{k+1}-\y_i^k\|^2\notag\\
=& f_i(\y_i^k)\hspace*{-0.05cm}+\hspace*{-0.05cm}\<\nabla f_i(\y_i^k),\x_i\hspace*{-0.05cm}-\hspace*{-0.05cm}\y_i^k\>\hspace*{-0.05cm}+\hspace*{-0.05cm}\<\nabla f_i(\y_i^k),\x_i^{k+1}\hspace*{-0.05cm}-\hspace*{-0.05cm}\x_i\>\hspace*{-0.05cm}+\hspace*{-0.05cm}\frac{L_i}{2}\|\x_i^{k+1}\hspace*{-0.05cm}-\hspace*{-0.05cm}\y_i^k\|^2\notag\\
\leq& f_i(\x_i)+\<\nabla f_i(\y_i^k),\x_i^{k+1}-\x_i\>+\frac{L_i}{2}\|\x_i^{k+1}-\y_i^k\|^2.\notag
\end{align} \end{subequations}
Adding the above two inequalities, we can have
\begin{subequations} \begin{align}
&F(\x^{k+1})-F(\x)\notag\\
\leq&\sum_{i=1}^2 \left[\<\A_i^T\overline\lambda_i^{k+1},\x_i-\x_i^{k+1}\>+\left(L_i+\frac{\beta\|\A_i\|_2^2 }{\theta^k}\right)\<\x_i^{k+1}-\y_i^k,\x_i-\y_i^k\>\right.\notag\\
&\left.-\left(\frac{L_i}{2}+\frac{\beta\|\A_i\|_2^2 }{\theta^k}\right)\|\x_i^{k+1}-\y_i^k\|^2\right].\notag
\end{align} \end{subequations}
Letting $\x_i=\x_i^k$ and $\x_i=\x_i^*$ respectively, we have
\begin{subequations} \begin{align}
&F(\x^{k+1})-F(\x^k)\notag\\
\leq&\sum_{i=1}^2 \left[\<\A_i^T\overline\lambda_i^{k+1},\x_i^k-\x_i^{k+1}\>+\left(L_i+\frac{\beta\|\A_i\|_2^2 }{\theta^k}\right)\<\x_i^{k+1}-\y_i^k,\x_i^k-\y_i^k\>\right.\notag\\
&\left.-\left(\frac{L_i}{2}+\frac{\beta\|\A_i\|_2^2 }{\theta^k}\right)\|\x_i^{k+1}-\y_i^k\|^2\right],\notag
\end{align} \end{subequations}
and
\begin{subequations} \begin{align}
&F(\x^{k+1})-F(\x^*)\notag\\
\leq&\sum_{i=1}^2 \left[\<\A_i^T\overline\lambda_i^{k+1},\x_i^*-\x_i^{k+1}\>+\left(L_i+\frac{\beta\|\A_i\|_2^2 }{\theta^k}\right)\<\x_i^{k+1}-\y_i^k,\x_i^*-\y_i^k\>\right.\notag\\
&\left.-\left(\frac{L_i}{2}+\frac{\beta\|\A_i\|_2^2 }{\theta^k}\right)\|\x_i^{k+1}-\y_i^k\|^2\right].\notag
\end{align} \end{subequations}
Multiplying the first inequality by $1-\theta^k$, multiplying the second by $\theta^k$ and adding them together, we have
\begin{subequations} \begin{align}
&F(\x^{k+1})-(1-\theta^k)F(\x^k)-\theta^kF(\x^*)\notag\\
\leq&\sum_{i=1}^2\left[\<\overline\lambda_i^{k+1},\theta^k\A_i\x_i^*+(1-\theta^k)\A_i\x_i^k-\A_i\x_i^{k+1}\>\right.\notag\\
&\left.+\left(L_i+\frac{\beta\|\A_i\|_2^2 }{\theta^k}\right)\<\x_i^{k+1}-\y_i^k,\theta^k\x_i^*+(1-\theta^k)\x_i^k-\y_i^k\>\right.\notag\\
&\left.-\left(\frac{L_i}{2}+\frac{\beta\|\A_i\|_2^2 }{\theta^k}\right)\|\x_i^{k+1}-\y_i^k\|^2\right].\notag
\end{align} \end{subequations}
Adding term $\<\lambda^*,\sum_{i=1}^2\A_i\x_i^{k+1}-(1-\theta^k)\sum_{i=1}^2\A_i\x_i^k-\theta^k\b\>$ to both sides, we can have
\begin{subequations} \begin{align}
&F(\x^{k+1})\hspace*{-0.07cm}-\hspace*{-0.07cm}F(\x^*)\hspace*{-0.07cm}+\hspace*{-0.07cm}\<\lambda^*,\A\x^{k+1}\hspace*{-0.07cm}-\hspace*{-0.07cm}\b\>\hspace*{-0.07cm}-\hspace*{-0.07cm}(1\hspace*{-0.07cm}-\hspace*{-0.07cm}\theta^k)\left( F(\x^k)\hspace*{-0.07cm}-\hspace*{-0.07cm}F(\x^*)\hspace*{-0.07cm}+\hspace*{-0.07cm}\<\lambda^*,\A\x^k\hspace*{-0.07cm}-\hspace*{-0.07cm}\b\> \right)\notag\\
=&F(\x^{k+1})-(1-\theta^k)F(\x^k)-\theta^kF(\x^*)\notag\\
&+\<\lambda^*,\sum_{i=1}^2\A_i\x_i^{k+1}-(1-\theta^k)\sum_{i=1}^2\A_i\x_i^k-\theta^k\b\>\notag\\
\leq&\sum_{i=1}^2\left[\<\overline\lambda_i^{k+1}-\lambda^*,\theta^k\A_i\x_i^*+(1-\theta^k)\A_i\x_i^k-\A_i\x_i^{k+1}\>\right.\notag\\
&\left.+\left(L_i+\frac{\beta\|\A_i\|_2^2 }{\theta^k}\right)\<\x_i^{k+1}-\y_i^k,\theta^k\x_i^*+(1-\theta^k)\x_i^k-\y_i^k\>\right.\notag\\
&\left.-\left(\frac{L_i}{2}+\frac{\beta\|\A_i\|_2^2 }{\theta^k}\right)\|\x_i^{k+1}-\y_i^k\|^2\right]\notag\\
=&\<\overline\lambda_1^{k+1}-\overline\lambda_2^{k+1},\theta^k\A_1\x_1^*+(1-\theta^k)\A_1\x_1^k-\A_1\x_1^{k+1}\>\notag\\
&+\sum_{i=1}^2\left[\<\overline\lambda_2^{k+1}-\lambda^*,\theta^k\A_i\x_i^*+(1-\theta^k)\A_i\x_i^k-\A_i\x_i^{k+1}\>\right.\notag\\
&\left.+\left(L_i+\frac{\beta\|\A_i\|_2^2 }{\theta^k}\right)\<\x_i^{k+1}-\y_i^k,\theta^k\x_i^*+(1-\theta^k)\x_i^k-\y_i^k\>\right.\notag\\
&\left.-\left(\frac{L_i}{2}+\frac{\beta\|\A_i\|_2^2 }{\theta^k}\right)\|\x_i^{k+1}-\y_i^k\|^2\right],\notag
\end{align} \end{subequations}
where we use $\sum_{i=1}^2\A_i\x_i^*=\b$. Let $\d_i^{k+1}=\frac{\x_i^{k+1}}{\theta^k}-\frac{1-\theta^k}{\theta^k}\x_i^k$ and $\d_i^k=\frac{\y_i^k}{\theta^k}-\frac{1-\theta^k}{\theta^k}\x_i^k,i=1,2$. Then we can have $\frac{\y_i^k}{\theta^k}-\frac{1-\theta^k}{\theta^k}\x_i^k=\frac{\x_i^k}{\theta^{k-1}}-\frac{1-\theta^{k-1}}{\theta^{k-1}}\x_i^{k-1}$, which leads to
\begin{subequations} \begin{align}
\y_i^k=\x_i^k+\frac{\theta^k(1-\theta^{k-1})}{\theta^{k-1}}(\x_i^k-\x_i^{k-1}).\notag
\end{align} \end{subequations}
which is (\ref{step1_LADMM_sn1}). From the definitions of $\overline\lambda_1^{k+1}$, $\overline\lambda_2^{k+1}$, $\d_i^{k+1}$ and $\d_i^k$, we can have
\begin{subequations} \begin{align}
&\<\overline\lambda_1^{k+1}-\overline\lambda_2^{k+1},\theta^k\A_1\x_1^*+(1-\theta^k)\A_1\x_1^k-\A_1\x_1^{k+1}\>\notag\\
=&\frac{\beta}{\theta^k}\<\A_1\y_1^k-\A_1\x_1^{k+1},\theta^k\A_1\x_1^*+(1-\theta^k)\A_1\x_1^k-\A_1\x_1^{k+1}\>\notag\\
=&\frac{\beta}{2\theta^k}\left[\|\theta^k\A_1\x_1^*\hspace*{-0.07cm}+\hspace*{-0.07cm}(1\hspace*{-0.07cm}-\hspace*{-0.07cm}\theta^k)\A_1\x_1^k\hspace*{-0.07cm}-\hspace*{-0.07cm}\A_1\x_1^{k+1}\|^2\hspace*{-0.07cm}-\hspace*{-0.07cm}\|\theta^k\A_1\x_1^*\hspace*{-0.07cm}+\hspace*{-0.07cm}(1\hspace*{-0.07cm}-\hspace*{-0.07cm}\theta^k)\A_1\x_1^k\hspace*{-0.07cm}-\hspace*{-0.07cm}\A_1\y_1^k\|^2 \right]\notag\\
&+\frac{\beta}{2\theta^k}\|\A_1\y_1^k-\A_1\x_1^{k+1}\|^2\notag\notag\\
\leq&\frac{\beta\theta^k}{2}\left[\|\A_1\d_1^{k+1}-\A_1\x_1^*\|^2-\|\A_1\d_1^k-\A_1\x_1^*\|^2 \right]+\frac{\beta\|\A_1\|_2^2}{2\theta^k}\|\y_1^k-\x_1^{k+1}\|^2,\notag
\end{align} \end{subequations}
and
\begin{subequations} \begin{align}
&\left(L_i+\frac{\beta\|\A_i\|_2^2 }{\theta^k}\right)\<\x_i^{k+1}-\y_i^k,\theta^k\x_i^*+(1-\theta^k)\x_i^k-\y_i^k\>\notag\\
=&\left(\frac{L_i}{2}\hspace*{-0.05cm}+\hspace*{-0.05cm}\frac{\beta\|\A_i\|_2^2 }{2\theta^k}\right)\left[\|\theta^k\x_i^*\hspace*{-0.05cm}+\hspace*{-0.05cm}(1\hspace*{-0.05cm}-\hspace*{-0.05cm}\theta^k)\x_i^k\hspace*{-0.05cm}-\hspace*{-0.05cm}\y_i^k\|^2\hspace*{-0.05cm}-\hspace*{-0.05cm}\|\theta^k\x_i^*\hspace*{-0.05cm}+\hspace*{-0.05cm}(1\hspace*{-0.05cm}-\hspace*{-0.05cm}\theta^k)\x_i^k\hspace*{-0.05cm}-\hspace*{-0.05cm}\x_i^{k+1}\|^2\right]\notag\\
&+\left(\frac{L_i}{2}+\frac{\beta\|\A_i\|_2^2 }{2\theta^k}\right)\|\x_i^{k+1}-\y_i^k\|^2.\notag\\
=&\frac{\theta^k\eta_i^k}{2}\left[\|\d_i^k-\x_i^*\|^2-\|\d_i^{k+1}-\x_i^*\|^2\right]+\left(\frac{L_i}{2}+\frac{\beta\|\A_i\|_2^2 }{2\theta^k}\right)\|\x_i^{k+1}-\y_i^k\|^2,\notag
\end{align} \end{subequations}
where $\eta_i^k=L_i\theta^k+\beta\|\A_i\|_2^2$. From Lemma \ref{lemma_lambda} we have
\begin{subequations} \begin{align}
&F(\x^{k+1})\hspace*{-0.07cm}-\hspace*{-0.07cm}F(\x^*)\hspace*{-0.07cm}+\hspace*{-0.07cm}\<\lambda^*,\A\x^{k+1}\hspace*{-0.07cm}-\hspace*{-0.07cm}\b\>\hspace*{-0.07cm}-\hspace*{-0.07cm}(1\hspace*{-0.07cm}-\hspace*{-0.07cm}\theta^k)\left( F(\x^k)\hspace*{-0.07cm}-\hspace*{-0.07cm}F(\x^*)\hspace*{-0.07cm}+\hspace*{-0.07cm}\<\lambda^*,\A\x^k\hspace*{-0.07cm}-\hspace*{-0.07cm}\b\> \right)\notag\\
\leq&\frac{\theta^k}{\beta}\<\overline\lambda_2^{k+1}-\lambda^*,\hat\lambda^k-\hat\lambda^{k+1}\>\notag\\
&+\frac{\beta\theta^k}{2}\left[\|\A_1\d_1^{k+1}-\A_1\x_1^*\|^2-\|\A_1\d_1^k-\A_1\x_1^*\|^2 \right]\notag\\
&+\theta^k\sum_{i=1}^2\frac{\eta_i^k}{2}\left[\|\d_i^k-\x_i^*\|^2-\|\d_i^{k+1}-\x_i^*\|^2\right]-\frac{\beta\|\A_2\|_2^2}{2\theta^k}\|\y_2^k-\x_2^{k+1}\|^2 \notag\\
=&\frac{\theta^k}{2\beta}\left(\|\hat\lambda^k-\lambda^*\|^2-\|\hat\lambda^{k+1}-\lambda^*\|^2-\|\overline\lambda_2^{k+1}-\hat\lambda^k\|^2+\|\hat\lambda^{k+1}-\overline\lambda_2^{k+1}\|^2\right)\notag\\
&+\frac{\beta\theta^k}{2}\left[\|\A_1\d_1^{k+1}-\A_1\x_1^*\|^2-\|\A_1\d_1^k-\A_1\x_1^*\|^2 \right]\notag\\
&+\theta^k\sum_{i=1}^2\frac{\eta_i^k}{2}\left[\|\d_i^k-\x_i^*\|^2-\|\d_i^{k+1}-\x_i^*\|^2\right]-\frac{\beta\|\A_2\|_2^2}{2\theta^k}\|\y_2^k-\x_2^{k+1}\|^2\notag\\
\leq&\frac{\theta^k}{2\beta}\left(\|\hat\lambda^k-\lambda^*\|^2-\|\hat\lambda^{k+1}-\lambda^*\|^2-\|\overline\lambda_2^{k+1}-\hat\lambda^k\|^2\right)\notag\\
&+\frac{\beta\theta^k}{2}\left[\|\A_1\d_1^{k+1}-\A_1\x_1^*\|^2-\|\A_1\d_1^k-\A_1\x_1^*\|^2 \right]\notag\\
&+\theta^k\sum_{i=1}^2\frac{\eta_i^k}{2}\left[\|\d_i^k-\x_i^*\|^2-\|\d_i^{k+1}-\x_i^*\|^2\right].\notag
\end{align} \end{subequations}
Dividing both sides by $\theta^k$ and using $\frac{1-\theta^k}{\theta^k}=\frac{1}{\theta^{k-1}}-\tau$, we have
\begin{subequations} \begin{align}
&\frac{1}{\theta^k}\left(F(\x^{k+1})\hspace*{-0.1cm}-\hspace*{-0.1cm}F(\x^*)\hspace*{-0.1cm}+\hspace*{-0.1cm}\<\lambda^*,\A\x^{k+1}\hspace*{-0.1cm}-\hspace*{-0.1cm}\b\>\right)\hspace*{-0.1cm}-\hspace*{-0.1cm}\frac{1}{\theta^{k-1}}\left( F(\x^k)\hspace*{-0.1cm}-\hspace*{-0.1cm}F(\x^*)\hspace*{-0.1cm}+\hspace*{-0.1cm}\<\lambda^*,\A\x^k\hspace*{-0.1cm}-\hspace*{-0.1cm}\b\> \right)\notag\\
&+\tau\left( F(\x^k)-F(\x^*)+\<\lambda^*,\A\x^k-\b\> \right)\notag\\
\leq&\frac{1}{2\beta}\left(\|\hat\lambda^k-\lambda^*\|^2-\|\hat\lambda^{k+1}-\lambda^*\|^2-\|\overline\lambda_2^{k+1}-\hat\lambda^k\|^2\right)\notag\\
&+\hspace*{-0.1cm}\left(\frac{\eta_1^k}{2}\|\d_1^k\hspace*{-0.1cm}-\hspace*{-0.1cm}\x_1^*\|^2\hspace*{-0.1cm}-\hspace*{-0.1cm}\frac{\beta}{2}\|\A_1\d_1^k\hspace*{-0.1cm}-\hspace*{-0.1cm}\A_1\x_1^*\|^2\hspace*{-0.1cm}-\hspace*{-0.1cm}\frac{\eta_1^k}{2}\|\d_1^{k+1}\hspace*{-0.1cm}-\hspace*{-0.1cm}\x_1^*\|^2\hspace*{-0.1cm}+\hspace*{-0.1cm}\frac{\beta}{2}\|\A_1\d_1^{k+1}\hspace*{-0.1cm}-\hspace*{-0.1cm}\A_1\x_1^*\|^2\right)\notag\\
&+\frac{\eta_2^k}{2}\left[\|\d_2^k-\x_2^*\|^2-\|\d_2^{k+1}-\x_2^*\|^2\right]\notag\\
\leq&\frac{1}{2\beta}\left(\|\hat\lambda^k-\lambda^*\|^2-\|\hat\lambda^{k+1}-\lambda^*\|^2-\|\overline\lambda_2^{k+1}-\hat\lambda^k\|^2\right)\notag\\
&+\hspace*{-0.15cm}\left(\hspace*{-0.1cm}\frac{\eta_1^k}{2}\hspace*{-0.03cm}\|\d_1^k\hspace*{-0.1cm}-\hspace*{-0.1cm}\x_1^*\|^2\hspace*{-0.1cm}-\hspace*{-0.1cm}\frac{\beta}{2}\hspace*{-0.03cm}\|\A_1\d_1^k\hspace*{-0.1cm}-\hspace*{-0.1cm}\A_1\x_1^*\|^2\hspace*{-0.1cm}-\hspace*{-0.1cm}\frac{\eta_1^{k+1}}{2}\hspace*{-0.03cm}\|\d_1^{k+1}\hspace*{-0.1cm}-\hspace*{-0.1cm}\x_1^*\|^2\hspace*{-0.1cm}+\hspace*{-0.1cm}\frac{\beta}{2}\hspace*{-0.03cm}\|\A_1\d_1^{k+1}\hspace*{-0.1cm}-\hspace*{-0.1cm}\A_1\x_1^*\|^2\hspace*{-0.15cm}\right)\notag\\
&+\frac{\eta_2^k}{2}\|\d_2^k-\x_2^*\|^2-\frac{\eta_2^{k+1}}{2}\|\d_2^{k+1}-\x_2^*\|^2,\notag
\end{align} \end{subequations}
where we use $\theta^{k+1}\leq \theta^k$ and $\eta_i^{k+1}\leq \eta_i^k$, which can be derived from $\frac{1}{\theta^{k+1}}-1=\frac{1}{\theta^{k}}-\tau$ and $0< \tau< 1$.$\hfill\blacksquare$
\end{proof}

A good property of Proposition \ref{serial_gc_general_convex_decent_lemma} is that we can sum the inequality over $k=0,\cdots,K$ and then bound $\frac{1}{\theta^K}(F(\x^{K+1})-F(\x^*)+\langle\lambda^*,\A\x^{K+1}-\b\rangle)$ by a constant, which leads to $F(\x^{K+1})-F(\x^*)+\langle\lambda^*,\A\x^{K+1}-\b\rangle\leq O(\theta^K)$. For the constraint functions, we have a similar result, which is described in the following proposition.
\begin{proposition}\label{constraint_proposition}
If the conditions in Proposition \ref{serial_gc_general_convex_decent_lemma} hold, then for Algorithm \ref{gc_serial_alg} we have
\begin{equation}
\left\|\sum_{k=0}^K \left(\frac{\A\x^{k+1}-\b}{\theta^k}-\frac{\A\x^{k}-\b}{\theta^{k-1}}+\tau\left(\A\x^{k}-\b\right)\right)\right\|\leq \frac{\sqrt{2\beta C}+\|\lambda^*-\hat\lambda^0\|}{\beta},\notag
\end{equation}
where $C=\frac{1}{2\beta}\|\lambda^0-\lambda^*\|^2+\frac{L_1+\beta\|\A_1\|_2^2}{2}\|\x_1^0-\x_1^*\|^2-\frac{\beta}{2}\|\A_1\x_1^0-\A_1\x_1^*\|^2+\frac{L_2+\beta\|\A_2\|_2^2}{2}\|\x_2^0$
$-\x_2^*\|^2$.
\end{proposition}
\begin{proof}
Summing the inequality in Proposition \ref{serial_gc_general_convex_decent_lemma} over $k=0,1,\cdots,K$, we have
\begin{eqnarray}
\begin{aligned}\label{sum_bound}
&\frac{1}{\theta^K}\left(F(\x^{K+1})-F(\x^*)+\<\lambda^*,\A\x^{K+1}-\b\>\right)\\
&+\sum_{k=1}^K \tau\left( F(\x^k)-F(\x^*)+\<\lambda^*,\A\x^k-\b\> \right)\leq C-\frac{1}{2\beta}\|\hat\lambda^{K+1}-\lambda^*\|^2,
\end{aligned}
\end{eqnarray}
where we use $\theta^0=1$, $0=\frac{1-\theta^0}{\theta^0}=\frac{1}{\theta^{-1}}-\tau$,
\begin{subequations} \begin{align}
\frac{\eta_1^{K+1}}{2}\|\d_1^{K+1}-\x_1^*\|^2-\frac{\beta}{2}\|\A_1\d_1^{K+1}-\A_1\x_1^*\|^2\geq 0,\notag
\end{align} \end{subequations}
and
\begin{subequations} \begin{align}
C\equiv&\frac{1}{2\beta}\|\lambda^0-\lambda^*\|^2+\left(\frac{L_1+\beta\|\A_1\|_2^2}{2}\|\x_1^0-\x_1^*\|^2-\frac{\beta}{2}\|\A_1\x_1^0-\A_1\x_1^*\|^2\right)\notag\\
&+\frac{L_2+\beta\|\A_2\|_2^2}{2}\|\x_2^0-\x_2^*\|^2\notag\\
=&\frac{1}{2\beta}\|\hat\lambda^0-\lambda^*\|^2+\left(\frac{\eta_1^0}{2}\|\d_1^0-\x_1^*\|^2-\frac{\beta}{2}\|\A_1\d_1^0-\A_1\x_1^*\|^2\right)+\frac{\eta_2^0}{2}\|\d_2^0-\x_2^*\|^2.\notag
\end{align} \end{subequations}
The last relation comes from $\d_i^0=\x_i^0$, $\hat\lambda^0=\lambda^0+\frac{\beta(1-\theta^0)}{\theta^0}\left( \sum_{i=1}^2\A_i\x_i^0-\b\right)=\lambda^0$ and $\eta_i^0=L_i\theta^0+\beta\|\A_i\|_2^2$.

Since $\{\x^*,\lambda^*\}$ is any KKT point, we have
\begin{subequations} \begin{align}
\x^*=\argmin_{\x} F(\x)+\<\lambda^*,\sum_{i=1}^2\A_i\x_i-\b\>.\notag
\end{align} \end{subequations}
So
\begin{equation}\label{function_bound}
F(\x^*)\hspace*{-0.05cm}=\hspace*{-0.05cm}F(\x^*)\hspace*{-0.05cm}+\hspace*{-0.05cm}\<\lambda^*,\sum_{i=1}^2\A_i\x_i^*\hspace*{-0.05cm}-\hspace*{-0.05cm}\b\>\leq F(\x)\hspace*{-0.05cm}+\hspace*{-0.05cm}\<\lambda^*,\sum_{i=1}^2\A_i\x_i\hspace*{-0.05cm}-\hspace*{-0.05cm}\b\>,\forall \x.
\end{equation}
Thus we have
\begin{subequations} \begin{align}
\frac{1}{2\beta}\|\hat\lambda^{K+1}-\lambda^*\|^2\leq C,\notag
\end{align} \end{subequations}
which leads to
\begin{subequations} \begin{align}
\|\hat\lambda^{K+1}-\hat\lambda^0\|\leq\|\hat\lambda^{K+1}-\lambda^*\|+\|\lambda^*-\hat\lambda^0\|\leq \sqrt{2\beta C}+\|\lambda^*-\hat\lambda^0\|.\notag
\end{align} \end{subequations}
From Lemma \ref{lemma_lambda}, we have
\begin{subequations} \begin{align}
&\left\|\sum_{k=0}^K \left[ \frac{\sum_{i=1}^2\A_i\x_i^{k+1}-\b}{\theta^k}-\frac{\sum_{i=1}^2\A_i\x_i^{k}-\b}{\theta^{k-1}}+\tau\left(\sum_{i=1}^2\A_i\x_i^{k}-\b\right)\right]\right\|\notag\\
\leq& \frac{\sqrt{2\beta C}+\|\lambda^*-\hat\lambda^0\|}{\beta}.\notag
\end{align} \end{subequations}
$\hfill\blacksquare$
\end{proof}

Both Propositions \ref{serial_gc_general_convex_decent_lemma} and \ref{constraint_proposition} have a similar form to (\ref{APG-descent}). Thus we have extended Nesterov's first acceleration scheme from unconstrained problems to constrained problems. Moreover, from Proposition \ref{constraint_proposition} we can see that Nesterov's acceleration scheme is critical to accelerate not only the decrease of the objective, but also the constraint error.

In Proposition \ref{constraint_proposition}, the summation lies inside the norm $\|\cdot\|$. Thus it is more difficult to bound $\left\|\frac{\A\x^{K+1}-\b}{\theta^K}\right\|$ than bounding $\frac{1}{\theta^K}(F(\x^{K+1})-F(\x^*)+\langle\lambda^*,\A\x^{K+1}-\b\rangle)$ from Propositon \ref{serial_gc_general_convex_decent_lemma}. We discover the following critical Lemma which can overcome this difficulty.
\begin{lemma}\label{bound_lemma}
Consider a sequence $\{\a^1,\a^2,\cdots\}$ of vectors, if $\{\a^k\}$ satisfies
\begin{subequations} \begin{align}
\left\|(1/\tau+K(1/\tau-1))\a^{K+1}+\sum_{k=1}^K\a^k\right\|\leq  c,\quad \forall K=0,1,2,\cdots.\notag
\end{align} \end{subequations}
where $1>\tau>0$. Then
$\|\sum_{k=1}^{K}\a^k\|< c$ for all $K=1,2,\cdots$.
\end{lemma}
\begin{proof}
For each $K\geq 0$, there exists $\c^{K+1}$ with every entry $\c_i^{K+1}\geq 0$ such that
\begin{subequations} \begin{align}
-\c_i^{K+1} \leq (1/\tau+K(1/\tau-1))\a_i^{K+1}+\sum_{k=1}^K\a_i^k \leq \c_i^{K+1},\notag
\end{align} \end{subequations}
and $\|\c^{K+1}\|=c$. Let $\s_i^K=\sum_{k=1}^K\a_i^k,\forall K\geq 1$ and $\s_i^0=0$, then
\begin{subequations} \begin{align}
\frac{-\c_i^{K+1}-\s_i^K}{1/\tau+K(1/\tau-1)}\leq \a_i^{K+1}\leq \frac{\c_i^{K+1}-\s_i^K}{1/\tau+K(1/\tau-1)},\forall K\geq 0,\notag
\end{align} \end{subequations}
where we use $1/\tau>1$ and $1/\tau+K(1/\tau-1)>0$. Thus, for all $K\geq 0$, we have
\begin{subequations} \begin{align}
&\s_i^{K+1}\notag\\
=&\a_i^{K+1}+\s_i^K\notag\\
\leq& \frac{\c_i^{K+1}-\s_i^K}{1/\tau+K(1/\tau-1)}+\s_i^K\notag\\
=&\frac{\c_i^{K+1}}{1/\tau+K(1/\tau-1)}+\frac{(K+1)(1/\tau-1)}{1/\tau+K(1/\tau-1)}\s_i^K\notag\\
\leq&\hspace*{-0.1cm}\frac{\c_i^{K+1}}{1/\tau+K(1/\tau-1)}\notag\\
&+\hspace*{-0.05cm}\frac{(K\hspace*{-0.05cm}+\hspace*{-0.05cm}1)(1/\tau\hspace*{-0.05cm}-\hspace*{-0.05cm}1)}{1/\tau\hspace*{-0.05cm}+\hspace*{-0.05cm}K(1/\tau\hspace*{-0.05cm}-\hspace*{-0.05cm}1)}\hspace*{-0.05cm}\left(\hspace*{-0.05cm} \frac{\c_i^{K}}{1/\tau\hspace*{-0.05cm}+\hspace*{-0.05cm}(K\hspace*{-0.05cm}-\hspace*{-0.05cm}1)(1/\tau\hspace*{-0.05cm}-\hspace*{-0.05cm}1)}\hspace*{-0.05cm}+\hspace*{-0.05cm}\frac{K(1/\tau\hspace*{-0.05cm}-\hspace*{-0.05cm}1)}{1/\tau\hspace*{-0.05cm}+\hspace*{-0.05cm}(K\hspace*{-0.05cm}-\hspace*{-0.05cm}1)(1/\tau\hspace*{-0.05cm}-\hspace*{-0.05cm}1)}\s_i^{K-1} \right)\notag\\
\leq&\frac{\c_i^{K+1}}{1/\tau+K(1/\tau-1)}+\frac{(K+1)(1/\tau-1)}{1/\tau+K(1/\tau-1)}\frac{\c_i^{K}}{1/\tau+(K-1)(1/\tau-1)}\notag\\
&+\frac{(K+1)(1/\tau-1)}{1/\tau+K(1/\tau-1)}\frac{K(1/\tau-1)}{1/\tau+(K-1)(1/\tau-1)}\left( \frac{\c_i^{K-1}}{1/\tau+(K-2)(1/\tau-1)}\right.\notag\\
&\left.+\frac{(K-1)(1/\tau-1)}{1/\tau+(K-2)(1/\tau-1)}\s_i^{K-2} \right)\notag\\
\leq&\frac{\c_i^{K+1}}{1/\tau+K(1/\tau-1)}\notag\\
&+\frac{(K+1)(1/\tau-1)}{1/\tau+K(1/\tau-1)}\frac{\c_i^{K}}{1/\tau+(K-1)(1/\tau-1)}\notag\\
&+\frac{(K+1)(1/\tau-1)}{1/\tau+K(1/\tau-1)}\frac{K(1/\tau-1)}{1/\tau+(K-1)(1/\tau-1)} \frac{\c_i^{K-1}}{1/\tau+(K-2)(1/\tau-1)}\notag\\
&+\cdots\notag\\
&+\hspace*{-0.05cm}\left(\hspace*{-0.05cm}\prod_{j=2}^{K+1}\frac{j(1/\tau\hspace*{-0.05cm}-\hspace*{-0.05cm}1)}{1/\tau\hspace*{-0.05cm}+\hspace*{-0.05cm}(j\hspace*{-0.05cm}-\hspace*{-0.05cm}1)(1/\tau\hspace*{-0.05cm}-\hspace*{-0.05cm}1)}\right)\hspace*{-0.05cm}\left( \hspace*{-0.05cm}\frac{\c_i^1}{1/\tau\hspace*{-0.05cm}+\hspace*{-0.05cm}0(1/\tau\hspace*{-0.05cm}-\hspace*{-0.05cm}1)}\hspace*{-0.05cm}+\hspace*{-0.05cm}\frac{1/\tau\hspace*{-0.05cm}-\hspace*{-0.05cm}1}{1/\tau\hspace*{-0.05cm}+\hspace*{-0.05cm}0(1/\tau\hspace*{-0.05cm}-\hspace*{-0.05cm}1)}\s_i^0 \right)\notag\\
=& \sum_{k=1}^{K+1}\left(\frac{\c_i^{k}}{1/\tau+(k-1)(1/\tau-1)}\prod_{j=k+1}^{K+1}\frac{j(1/\tau-1)}{1/\tau+(j-1)(1/\tau-1)}\right),\notag
\end{align} \end{subequations}
where we set $\prod_{j=K+2}^{K+1}\frac{j(1/\tau-1)}{1/\tau+(j-1)(1/\tau-1)}=1$. Define
\begin{subequations} \begin{align}
r^k\hspace*{-0.05cm}=\hspace*{-0.05cm}\frac{1}{1/\tau\hspace*{-0.05cm}+\hspace*{-0.05cm}(k\hspace*{-0.05cm}-\hspace*{-0.05cm}1)(1/\tau\hspace*{-0.05cm}-\hspace*{-0.05cm}1)}\prod_{j=k+1}^{K+1}\frac{j(1/\tau\hspace*{-0.05cm}-\hspace*{-0.05cm}1)}{1/\tau\hspace*{-0.05cm}+\hspace*{-0.05cm}(j\hspace*{-0.05cm}-\hspace*{-0.05cm}1)(1/\tau\hspace*{-0.05cm}-\hspace*{-0.05cm}1)},\forall k=1,2,\cdots,K+1.\notag
\end{align} \end{subequations}
Then we have $r^k>0$ and $\s_i^{K+1}\leq \sum_{k=1}^{K+1}r^k\c_i^k$. Similarly, we also have $\s_i^{K+1}\geq-\sum_{k=1}^{K+1}r^k\c_i^k$. Thus
\begin{subequations} \begin{align}
|\s_i^{K+1}|\leq \sum_{k=1}^{K+1}r^k\c_i^k.\notag
\end{align} \end{subequations}
Define
\begin{subequations} \begin{align}
R^{K+1}=\sum_{k=1}^{K+1}\frac{1}{1/\tau+(k-1)(1/\tau-1)}\prod_{j=k+1}^{K+1}\frac{j(1/\tau-1)}{1/\tau+(j-1)(1/\tau-1)},\notag
\end{align} \end{subequations}
\begin{subequations} \begin{align}
R^{K}=\sum_{k=1}^{K}\frac{1}{1/\tau+(k-1)(1/\tau-1)}\prod_{j=k+1}^{K}\frac{j(1/\tau-1)}{1/\tau+(j-1)(1/\tau-1)},\notag
\end{align} \end{subequations}
and
\begin{subequations} \begin{align}
R^{1}=\sum_{k=1}^{1}\frac{1}{1/\tau+(k-1)(1/\tau-1)}\prod_{j=k+1}^{1}\frac{j(1/\tau-1)}{1/\tau+(j-1)(1/\tau-1)}=\tau.\notag
\end{align} \end{subequations}
Then we have
\begin{subequations} \begin{align}
&R^{K+1}\notag\\
=&\frac{1}{1/\tau+K(1/\tau-1)}+\sum_{k=1}^{K}\frac{1}{1/\tau+(k-1)(1/\tau-1)}\prod_{j=k+1}^{K+1}\frac{j(1/\tau-1)}{1/\tau+(j-1)(1/\tau-1)}\notag\\
=&\frac{1}{1/\tau+K(1/\tau-1)}\notag\\
&+\frac{(K+1)(1/\tau-1)}{1/\tau+K(1/\tau-1)}\sum_{k=1}^K\frac{1}{1/\tau+(k-1)(1/\tau-1)}\prod_{j=k+1}^K\frac{j(1/\tau-1)}{1/\tau+(j-1)(1/\tau-1)}\notag\\
=&\frac{1}{1/\tau+K(1/\tau-1)}+\frac{(K+1)(1/\tau-1)}{1/\tau+K(1/\tau-1)}R^K.\notag
\end{align} \end{subequations}

Next, we prove $R^K< 1,\forall K\geq 1$ by induction. It can be easily checked that $R^1=\tau< 1$. Assume that $R^K<1$ holds, then
\begin{subequations} \begin{align}
R^{K+1}<\frac{1}{1/\tau+K(1/\tau-1)}+\frac{(K+1)(1/\tau-1)}{1/\tau+K(1/\tau-1)}=1.\notag
\end{align} \end{subequations}
So by induction we can have $R^K< 1,\forall K\geq 1$.

So for any $K\geq 0$, we have
\begin{subequations} \begin{align}
(\s_i^{K+1})^2\leq \left(\sum_{k=1}^{K+1}r^k\right)^2\left(\frac{\sum_{k=1}^{K+1}r^k\c_i^k}{\sum_{k=1}^{K+1}r^k}\right)^2\leq\left(\sum_{k=1}^{K+1}r^k\right)^2\frac{\sum_{k=1}^{K+1}r^k(\c_i^k)^2}{\sum_{k=1}^{K+1}r^k}<\sum_{k=1}^{K+1}r^k(\c_i^k)^2,\notag
\end{align} \end{subequations}
where we use $\sum_{k=1}^{K+1}r^k=R^{K+1}<1$ and the Jensen inequality for $x^2$. So we have
\begin{equation}
\|\S^{K+1}\|^2=\sum_{i}(\s_i^{K+1})^2<\sum_{k=1}^{K+1}r^k\sum_{i}(\c_i^k)^2=\sum_{k=1}^{K+1}r^k c^2<c^2,\notag
\end{equation}
where we use $\|\c^k\|=c,\forall k\geq 1$. So $\|\sum_{k=1}^{K+1}\a^k\|=\|\S^{K+1}\|< c,\forall K\geq 0$.$\hfill\blacksquare$
\end{proof}

Based on Propositions \ref{serial_gc_general_convex_decent_lemma} and \ref{constraint_proposition}, we can have the $O(1/K)$ nonergodic convergence rate in Theorem \ref{serial_general_covnex_theorem}.
\begin{theorem}\label{serial_general_covnex_theorem}
If the conditions in Proposition \ref{serial_gc_general_convex_decent_lemma} hold, then for Algorithm \ref{gc_serial_alg} we have
\begin{subequations} \begin{align}
-\frac{2\tau C_1\|\lambda^*\|}{1+K(1-\tau)}\leq F(\x^{K+1})-F(\x^*)\leq \frac{C+2\tau C_1\|\lambda^*\|}{1+K(1-\tau)},\notag
\end{align} \end{subequations}
and
\begin{subequations} \begin{align}
\left\|\A\x^{K+1}-\b\right\|\leq \frac{2\tau C_1}{1+K(1-\tau)},\notag
\end{align} \end{subequations}
where $C_1=\frac{\sqrt{2\beta C}+\|\lambda^*-\lambda^0\|}{\tau\beta}$ and $C$ is defined in Proposition \ref{constraint_proposition}.
\end{theorem}
\begin{proof}
From (\ref{sum_bound}), (\ref{function_bound}) and Proposition \ref{constraint_proposition}, we can have
\begin{subequations} \begin{align}
F(\x^{K+1})-F(\x^*)+\<\lambda^*,\sum_{i=1}^2\A_i\x_i^{K+1}-\b\>\leq C\theta^K,\notag\
\end{align} \end{subequations}
and
\begin{subequations} \begin{align}
&\frac{\sqrt{2\beta C}+\|\lambda^*-\hat\lambda^0\|}{\beta}\notag\\
\geq&\left\|\sum_{k=0}^K \left\{ \frac{\sum_{i=1}^2\A_i\x_i^{k+1}-\b}{\theta^k}-\frac{\sum_{i=1}^2\A_i\x_i^{k}-\b}{\theta^{k-1}}+\tau\left(\sum_{i=1}^2\A_i\x_i^{k}-\b\right)\right\}\right\|\notag\\
=&\left\|\frac{\sum_{i=1}^2\A_i\x_i^{K+1}-\b}{\theta^K}-\frac{\sum_{i=1}^2\A_i\x_i^{0}-\b}{\theta^{-1}}+\sum_{k=0}^K\tau\left(\sum_{i=1}^2\A_i\x_i^{k}-\b\right)\right\|\notag\\
=&\left\|\frac{\sum_{i=1}^2\A_i\x_i^{K+1}-\b}{\theta^K}+\sum_{k=1}^K\tau\left(\sum_{i=1}^2\A_i\x_i^{k}-\b\right)\right\|,\forall K=0,1,2,\cdots\notag.
\end{align} \end{subequations}
where we use $\frac{1}{\theta^{-1}}-\tau=\frac{1-\theta^0}{\theta^0}=0$. Since $\frac{1}{\theta^k}=\frac{1}{\theta^{k-1}}+1-\tau= \frac{1}{\theta^0}+k(1-\tau)$, we have $\theta^k= \frac{1}{\frac{1}{\theta^0}+k(1-\tau)}=\frac{1}{1+k(1-\tau)}$. For simplicity, let $\a^k=\sum_{i=1}^2\A_i\x_i^{k}-\b$. Then we can have
\begin{subequations} \begin{align}
\left\|(1/\tau\hspace*{-0.05cm}+\hspace*{-0.05cm}K(1/\tau\hspace*{-0.05cm}-\hspace*{-0.05cm}1))\a^{K+1}\hspace*{-0.05cm}+\hspace*{-0.05cm}\sum_{k=1}^K\a^k\right\|\hspace*{-0.05cm}\leq\hspace*{-0.05cm}\frac{\sqrt{2\beta C}\hspace*{-0.05cm}+\hspace*{-0.05cm}\|\lambda^*\hspace*{-0.05cm}-\hspace*{-0.05cm}\hat\lambda^0\|}{\tau\beta}\equiv C_1, \forall K=0,1,\cdots.\notag
\end{align} \end{subequations}
From Lemma \ref{bound_lemma} we have $\|\sum_{k=1}^K\a^k\|\leq C_1,\forall K=1,2,\cdots$. So $\|\a^{K+1}\|\leq \frac{2C_1}{1/\tau+K(1/\tau-1)}$, $\forall K=1,2,\cdots$. Moreover, $\|\a^1\|\leq \tau C_1\leq \frac{2C_1}{1/\tau+0(1/\tau-1)}$. So
\begin{subequations} \begin{align}
\left\|\sum_{i=1}^2\A_i\x_i^{K+1}-\b\right\|\leq \frac{2\tau C_1}{1+K(1-\tau)},\forall K=0,1,\cdots,\notag
\end{align} \end{subequations}
Thus we can have
\begin{subequations} \begin{align}
F(\x^{K+1})-F(\x^*)\leq& C\theta^K+\|\lambda^*\|\left\|\sum_{i=1}^2\A_i\x_i^{K+1}-\b\right\|\notag\\
\leq& \frac{C}{1+K(1-\tau)}+\frac{2\tau C_1\|\lambda^*\|}{1+K(1-\tau)},\notag
\end{align} \end{subequations}
and
\begin{subequations} \begin{align}
F(\x^{K+1})-F(\x^*)\geq -\|\lambda^*\|\left\|\sum_{i=1}^2\A_i\x_i^{K+1}-\b\right\|\geq -\frac{2\tau C_1\|\lambda^*\|}{1+K(1-\tau)},\notag
\end{align} \end{subequations}
which is derived from (\ref{function_bound}).$\hfill\blacksquare$
\end{proof}

From Theorem \ref{serial_general_covnex_theorem} we can see that the $O(1/K)$ nonergodic convergence rate exists only if $\tau<1$. In fact, only when $\tau<1$, $\theta^k=\frac{1}{1+k(1-\tau)}$ is in the order of $O(1/k)$ and Nesterov's acceleration scheme is effective. As discussed in Section \ref{equi_sec}, ALADMM-NE reduces to the traditional LADMM when $\tau=1$.

\subsection{Tips on the Choice of the Algorithms}\label{tip_sec}
In applications where the practical performance of (L)ADMM coincides with its theoretical convergence rate, it is guaranteed that ALADMM-NE practically outperforms (L)ADMM. However, in the cases where (L)ADMM converges much faster than its theoretical rate, e.g., in applications of Robust PCA \cite{lin_alm} that (L)ADMM almost linearly converges, we empirically observe that the superiority of ALADMM-NE and the accelerated ADMM in \cite{ouyang-2015-fastLADM} is not obvious. In fact, due to the special setting of $\theta^k$ which dependents on $k$, ALADMM-NE and the method in \cite{ouyang-2015-fastLADM} have exactly the $O(1/K)$ convergence rate measured at $\{\x_1^K,\x_2^K\}$ even for the strongly convex problems. So in practice, we suggest that when the problem is complex and does not satisfy the linear convergence conditions \cite{Deng-2012-linearADM,Luo-2012-LinearADM,Giselsson-2016,han-2016,boley-2013}, ALADMM-NE and the accelerated ADMM in \cite{ouyang-2015-fastLADM} are better choices than the traditional (L)ADMM. When sparseness or low-rankness is required, ALADMM-NE is better than the accelerated ADMM in \cite{ouyang-2015-fastLADM}.

Donoghue and Cand\`{e}s \cite{Donoghue-2015-NesRestart} proposed a restart strategy for Nesterov's first acceleration scheme when minimizing the unconstrained problems, in which the algorithm is restarted after some iterations by setting $\theta^{k+1}=1$ and $\y^{k+1}=\x^{k+1}$. Then the linear convergence is guaranteed even for the sublinear setting of $\theta^k$ \cite{Necoara-2016}. A similar technique is discussed for Nesterov's second scheme in \cite{li-2017}. So we can apply the restart scheme for the accelerated ADMM in \cite{ouyang-2015-fastLADM} and ALADMM-NE. The latter is described in Algorithm \ref{gc_serial_algr}. We restart ALADMM-NE as long as $\|\A\x^{k+1}-\b\|$ increases.  We set $\theta^{k+1}=\theta^k=1$ in the if-clause to make $\y^{k+1}=\x^{k+1}$ when the algorithm is restarted. We use the criterion $\theta^{k+1}<\epsilon$ to prevent frequent restart and only restart when $\theta^k$ becomes small.

\begin{algorithm}[tb]
   \caption{Accelerated LADMM with NonErgodic convergence rate and Restart(ALADMM-NER)}
   \label{gc_serial_algr}
\begin{algorithmic}
   \STATE Initialize $\lambda^0$, $\x_i^0=\x_i^{-1},i=1,2$, $1>\tau> 0.5$, $\beta>0$, $\theta^0=1$, $1>\epsilon>0$, $\theta^{-1}=1/\tau$.
   \FOR{$k=0,1,2,\cdots$}
   \STATE Update $\y_i^k,i=1,2$ using (\ref{step1_LADMM_sn1}),\\
   \STATE Update $\x_1^{k+1}$ and $\x_2^{k+1}$ serially using (\ref{step2_LADMM_n1}) and (\ref{step3_LADMM_n1}),\\
   \STATE Update $\hat\lambda^{k+1}$ using (\ref{step4_LADMM_sn1}),\\
   $\theta^{k+1}=\frac{1}{1-\tau+\frac{1}{\theta^k}}$.
   \IF{$\|\sum_{i=1}^2\A_i\x_i^{k+1}-\b\|\geq \|\sum_{i=1}^2\A_i\x_i^k-\b\|$ and $\theta^{k+1}<\epsilon$}
   \STATE $\theta^{k+1}=1$, $\theta^k=1$\\
   \ENDIF
   \ENDFOR
\end{algorithmic}
\end{algorithm}

\section{Tightness of the $o(1/\sqrt{K})$ Nonergodic Rate for the Traditional ADMM}\label{ADMM-tight}

In this section we show that the $o\left(\frac{1}{\sqrt{K}}\right)$ rate is tight for ADMM, at least for the constraint,. We study a special problem \cite{Bauschke-2014-DRslow,Davis-2014-DR}, on which the Alternating Projection Method (APM) and DR splitting perform slowly. They converge arbitrarily slowly on the measure of $\|\x^k-\x^*\|$ and converge with the tight $o\left(\frac{1}{\sqrt{k}}\right)$ rate on the measure of $f(\x^k)-f(\x^*)$. The discussion in this section also suits for LADMM and the accelerated ADMM in \cite{ouyang-2015-fastLADM} (measured at $(\z_1^k,\z_2^k)$) since they are equivalent to ADMM on this special problem.

Let $\vartheta_i$ be a sequence of angles in $(0,\pi/2)$ with $\cos(\vartheta_i)=\sqrt{\frac{i}{i+1}}$. Let $\e_0=(1,0)$, $\e_{\pi/2}=(0,1)$ and $\e_{\vartheta_i}=\cos(\vartheta_i)\e_0+\sin(\vartheta_i)\e_{\pi/2}$. Define two lines $U=\mbox{span}\{\e_0\}$ and $V_i=\mbox{span}\{\e_{\vartheta_i}\}$, then $U \bigcap V_i=\{0\}$. Consider the Hilbert space $\mathbf H=\mathbf R^2\bigoplus\mathbf R^2\bigoplus\cdots$ and define
\begin{subequations} \begin{align}
&\U=\mathbf R\cdot \e_0\times \mathbf R\cdot \e_0\times\cdots,\notag\\
&\V=\mathbf R\cdot \e_{\vartheta_0}\times \mathbf R\cdot \e_{\vartheta_1}\times\cdots.\notag
\end{align} \end{subequations}
We consider problem
\begin{equation}\label{problem-hard}
\min_{\x} f(\x)=h(\x)+g(\x),
\end{equation}
where $h(\x)=\mbox{I}_{\U}(\x)$ is the indicator function of $\U$, $g(\x)=\frac{\beta}{\sqrt{2a-1}} \mbox{d}_{\V}(\x)$, $\mbox{d}_{\V}(\x)=\min_{\v\in \V}\|\x-\v\|$ and $a$ can be any constant satisfying $a>0.5$. This problem can be solved by ADMM and ALADMM-NE by transforming it to
\begin{equation}\label{problem-ADM}
\min_{\x,\z} h(\x)+g(\z)\qquad s.t. \quad\z-\x=0.
\end{equation}

Proposition \ref{theorem-bound} says that the $o\left(\frac{1}{\sqrt{K}}\right)$ rate is tight for ADMM. This means that the slow $o\left(\frac{1}{\sqrt{K}}\right)$ nonergodic convergence rate of ADMM is not due to the weakness of the proof, but that of ADMM itself. It is difficult to establish the lower complexity bound of $|h(\x^k)+g(\z^k)-h(\x^*)-g(\z^*)|$, so we only measure $f(\x^k)-f(\x^*)$ for simplicity. It should be noted that Proposition \ref{theorem-bound} is ADMM specified and it does not suit for ALADMM-NE. As comparison, we can establish $\|\z^{k+1}-\x^{k+1}\|\leq O(1/k)$ and $f(\x^k)-f(\x^*)\leq O(1/k)$ for ALADMM-NE, which establishes the superiority of ALADMM-NE with theoretical guarantee\footnote{ALADMM-NE can be applied to Hilbert spaces. Since $g(\z)$ is continuous \cite{Davis-2014-DR}, we have $f(\x^k)-f(\x^*)\leq |h(\x^k)+g(\z^k)-h(\x^*)-g(\z^*)|+|g(\z^k)-g(\x^k)|\leq O(1/k)+O(L/k)=O(1/k)$.}. We list the comparisons in Table \ref{table-hard}.
\begin{proposition}\label{theorem-bound}
Let $\x^0=\left(\left[
  \begin{array}{cc}
    \frac{1}{(i+1)^a} \\
    0
  \end{array}\notag
  \right]\right)_{i\geq 1}$, $\lambda^0=0$, $a>0.5$, then for ADMM with iterations (\ref{adm1})-(\ref{adm3}) we have
$\|\z^{k+1}-\x^{k+1}\|\geq \Omega\left(\frac{1}{(k+2)^a}\right)$ and $f(\x^k)-f(\x^*)\geq\Omega\left(\frac{\beta}{\sqrt{2a-1}(k+1)^a}\right)$.
\end{proposition}

In Proposition \ref{theorem-bound} we specialize the initialization of $\x^0$ and $\lambda^0$, where $\|\x^0-\x^*\|$ is bounded and independent on $k$. This is a standard trick in the analysis of lower bound. Proposition \ref{theorem-bound}
can be proved using the same proof framework in \cite{Davis-2014-DR}, so we  omit the details.

\begin{table}[t]
\caption{Theoretical complexity comparisons among ALADMM-NE, ADMM, DR and APM on problem (\ref{problem-hard}). $a$ is any constant satisfying $a>0.5$.\label{table-hard}}
\begin{center}
\begin{tabular}{c||c}\hline
  & Theoretical complexity bound
\\\hline
APM & $f(\x^k)-f(\x^*)\geq\Omega\left(\frac{1}{k^a}\right)$ \\
\hline
DR  & $f(\x^k)-f(\x^*)\geq\Omega\left(\frac{1}{k^a}\right)$ \\
\hline
ADMM & $f(\x^k)-f(\x^*)\geq\Omega\left(\frac{1}{k^a}\right)$, $\|\z^k-\x^k\|\geq\Omega\left(\frac{1}{k^a}\right)$ \\
\hline
ALADMM-NE & $f(\x^k)-f(\x^*)\leq O(1/k)$, $\|\z^k-\x^k\|\leq O(1/k)$\\
\hline
\end{tabular}
\end{center}
\end{table}

One may think that the increasing penalty $\frac{\beta}{\theta^k}$ in ALADMM-NE is the deciding factor of the improved convergence rate. However, this is incorrect. Empirically, large penalty speeds up the decrease of the constraint error in ADMM \cite{Lin-2015-LADMPSAP}. But this is not guaranteed in theory. In fact, From Proposition \ref{theorem-bound} we can see that the constraint error is independent of $\beta$, which means that the decrease of the constraint error cannot be faster than $o\left(\frac{1}{\sqrt{K}}\right)$ no matter how large $\beta$ is. There are two reasons for this result: 1. It is equivalent to minimizing the sum of two indicator functions when using ADMM to solve problem (\ref{problem-ADM}) and $\beta$ has no influence on the projection operation; 2. $\x$ and $\z$ are updated serially, not parallel. Thus although the gradually increasing penalty in ALADMM-NE plays an important role to cooperate with Nesterov's acceleration scheme, Nesterov's scheme is indeed the critical factor to improve the convergence rate in theory. Large penalty cannot improve the convergence rate of ADMM even for the constraint.

\section{Lower Complexity Bound}\label{lb_section}

Recently, Woodworth and Srebro \cite{Blake-ADMbound-2016} established the $O(1/K)$ lower complexity bound of the stochastic gradient methods for optimizing the finite sum problem: $\min_{\x} \frac{1}{m}\sum_{i=1}^m f_i(\x)$, where each $f_i$ is nonsmooth and non-strongly convex. In this section we use Woodworth and Srebro's result to analyze the general splitting scheme, and then extend it to the general ADMM type methods, which deal with the additional linear constraint.

\subsection{Splitting Scheme}

We consider the following problem:
\begin{subequations} \begin{align}
\min_{\x\in \mathbb X} F_1(\x)+F_2(\x).\notag
\end{align} \end{subequations}

We call a method belonging to the general splitting scheme if it has the form of
\begin{eqnarray}
\begin{aligned}\label{splitting_scheme}
&\mbox{Generate } \z_1^t\mbox{ based on }\{\x_1^{1:t},\x_2^{1:t},\z_1^{1:t-1},\z_2^{1:t},F_1(\z_1^{1:t-1}),F_2(\z_2^{1:t})\},\\
&\x_1^{t+1}=\mbox{Prox}_{F_1/\beta^t}(\z_1^t),\\
&\mbox{Generate } \z_2^{t+1}\mbox{ based on }\{\x_1^{1:t+1},\x_2^{1:t},\z_1^{1:t},\z_2^{1:t},F_1(\z_1^{1:t}),F_2(\z_2^{1:t})\}, \\
&\x_2^{t+1}=\mbox{Prox}_{F_2/\beta^t}(\z_2^{t+1}),
\end{aligned}
\end{eqnarray}
at the $t$-th iteration and $\beta^t$ is arbitrary. We denote $\x^{1:t}=\{\x^1,\cdots,\x^t\}$ and $F(\x^{1:t})=\{F(\x^1),\cdots,F(\x^t)\}$ for simplicity. In this general scheme, two proximal subproblems are solved alternatively and $\{\z_1^k,\z_2^{k+1}\}$ can be generated in any way, e.g., $\z_1^t\in\mbox{Span}\{\x_1^{1:t},\x_2^{1:t},\z_1^{1:t-1},\z_2^{1:t}\}$ and $\z_2^{t+1}\in\mbox{Span}\{\x_1^{1:t+1},\x_2^{1:t},\z_1^{1:t},\z_2^{1:t}\}$. The algorithm belonging to this scheme accesses the objectives $F_1$ and $F_2$ only through the oracle of $(\mbox{Prox}_{F_i/\beta^t}(\x),F_i(\x),i=1,2)$. It generates the next iterates of $\{\z_1^{t+1},\z_2^{t+2}\}$ based on the previous responses of the oracle. This general splitting scheme includes many famous splitting algorithms, such as DR splitting, which consists of the following steps:
\begin{subequations} \begin{align}
&\x_1^{t+1}=\mbox{Prox}_{F_1/\beta}(\z^t),\notag\\
&\x_2^{t+1}=\mbox{Prox}_{F_2/\beta}(2\x_1^{t+1}-\z_1^t),\notag\\
&\z^{t+1}=\z^t-\x_1^{t+1}+\x_2^{t+1}.\notag
\end{align} \end{subequations}

For this general splitting scheme, we can have the $O\left(1/K\right)$ lower bound, which is described in the following proposition. Note that we do not aim to construct a counterexample such that for all algorithms satisfying (\ref{splitting_scheme}), they converge slowly. Instead, for any algorithm satisfying (\ref{splitting_scheme}), we want to construct a counterexample such that it converges slowly. The counterexample is not algorithm independent.
\begin{proposition}\label{splitting theorem}
For any algorithm belonging to the general splitting scheme (\ref{splitting_scheme}), there exist convex and $L$-Lipschitz continuous functions $F_1$ and $F_2$ defined over $\mathbb X=\{\x\in \mathbb{R}^{6k+2}:\|\x\|\leq B\}$, such that
\begin{subequations} \begin{align}
F_1(\hat\x^k)+F_2(\hat\x^k)\geq \frac{LB}{8(k+1)},\notag
\end{align} \end{subequations}
where $\hat\x^k=\sum_{i=1}^k\alpha_1^i\x_1^i+\sum_{i=1}^k\alpha_2^i\x_2^i$, $\forall\alpha_1^i$ and $\forall\alpha_2^i$, $i=1,\cdots,k$.
\end{proposition}
Proposition \ref{splitting theorem} can be proved using the same analysis framework in \cite{Blake-ADMbound-2016}. We give the proof sketch for the reader's convenience. For the detailed analysis, please see \cite{Blake-ADMbound-2016}.

\noindent\textbf{Proof Sketch:} For any algorithm belonging to the splitting scheme (\ref{splitting_scheme}), we want to construct a hard function for witch the algorithm converges slowly. For simplicity, we let $L=1$ and $B=1$. Initialize $\z_2^1$, $\x_1^1$, $\x_2^1$, $\v^1$, $\v^0$ and $F_1^1=\frac{1}{\sqrt{2}}|b-\<\x,\v^0\>|+\frac{1}{4\sqrt{k}}|\<\x,\v^0\>-\<\x,\v^1\>|$ such that $\|\v^1\|=1$, $\|\v^0\|=1$ and $\v^1\bot\v^0$. We use an adversary strategy to construct the hard function, i.e., at the $t$-th iteration the algorithm quires the oracle with $(\z_1^t,\z_2^{t+1},\beta^t)$ and an adversary responses with an answer of $(\mbox{Prox}_{F_1^t/\beta^t}(\z_1^t),\mbox{Prox}_{F_2^t/\beta^t}(\z_2^{t+1}),F_1^t(\z_1^t),F_2^t(\z_2^{t+1}))$. The algorithm accesses the problem only through the oracle and it makes the decisions based on the previous responses of the oracle. The adversary constructs the hard function gradually based on the previous queries of the algorithm. Specifically, at the $t$-th iteration with $t=1,\cdots,k$, we perform the following steps:
    \begin{eqnarray}
    \begin{split}\label{splitting_scheme1}
        &\mbox{Algorithm: }\\
        &\hspace*{0.8cm}\mbox{Generate } \z_1^t\mbox{ based on }\{\x_1^{1:t},\x_2^{1:t},\z_1^{1:t-1},\z_2^{1:t},F_1(\z_1^{1:t-1}),F_2(\z_2^{1:t})\}, \\
        &\hspace*{0.8cm}\x_1^{t+1}=\mbox{Prox}_{F_1^t/\beta^t}(\z_1^t),\quad F_1(\z_1^t)=F_1^t(\z_1^t),\\
        &\mbox{Adversary: }\\
        &\hspace*{0.8cm}\mbox{Construct }\v^{2t}\bot \{\v^{0:2t-1},\z_1^{1:t},\z_2^{1:t},\x_1^{1:t},\x_2^{1:t}\}\mbox{ such that }\|\v^{2t}\|=1,\\
        &\hspace*{0.8cm}\mbox{Construct }F_2^t=\frac{1}{\sqrt{2}}|b-\<\x,\v^0\>|+\frac{1}{4\sqrt{k}}\sum_{r=1}^{t}|\<\x,\v^{2r-1}\>-\<\x,\v^{2r}\>|,\hspace*{-1cm}\\
        &\mbox{Algorithm: }\\
        &\hspace*{0.8cm}\mbox{Generate } \z_2^{t+1}\mbox{ based on }\{\x_1^{1:t+1},\x_2^{1:t},\z_1^{1:t},\z_2^{1:t},F_1(\z_1^{1:t}),F_2(\z_2^{1:t})\},\\
        &\hspace*{0.8cm}\x_2^{t+1}=\mbox{Prox}_{F_2^t/\beta^t}(\z_2^{t+1}),\quad F_2(\z_2^{t+1})=F_2^t(\z_2^{t+1}),\\
        &\mbox{Adversary: }\\
        &\hspace*{0.8cm}\mbox{Construct }\v^{2t+1}\bot\{\v^{0:2t},\z_1^{1:t},\z_2^{1:t+1}\}\mbox{ such that }\|\v^{2t+1}\|\hspace*{-0.1cm}=\hspace*{-0.1cm}1,\\
        &\hspace*{0.8cm}\mbox{Construct }F_1^{t+1}\hspace*{-0.05cm}=\hspace*{-0.05cm}\frac{1}{\sqrt{2}}|b\hspace*{-0.05cm}-\hspace*{-0.05cm}\<\x,\v^0\>\hspace*{-0.05cm}|\hspace*{-0.05cm}+\hspace*{-0.05cm}\frac{1}{4\sqrt{k}}\sum_{r=1}^{t+1}\hspace*{-0.05cm}|\hspace*{-0.05cm}\<\x,\v^{2r-2}\>\hspace*{-0.05cm}-\hspace*{-0.05cm}\<\x,\v^{2r-1}\>\hspace*{-0.05cm}|,\hspace*{-1cm}
    \end{split}
    \end{eqnarray}
where $F_1^t$ and $F_2^t$ are adaptive of the history iterates, i.e., we construct $F_1^t$ and $F_2^t$ based on the history iterates and they are different from each other at different iterations. However, due to the orthogonality between $\v$ and $\z$, we can prove the following relations
    \begin{eqnarray}
        &&\mbox{Prox}_{F_1^t/\beta^t}(\z_1^t)=\mbox{Prox}_{F_1^{k}/\beta^t}(\z_1^t),\quad F_1^{t}(\z_1^t)=F_1^{k}(\z_1^t),\quad\forall t\leq k,\notag\\
        &&\mbox{Prox}_{F_2^t/\beta^t}(\z_2^{t+1})=\mbox{Prox}_{F_2^{k}/\beta^t}(\z_2^{t+1}),\quad F_2^t(\z_2^{t+1})=F_2^k(\z_2^{t+1}),\quad\forall t\leq k.\notag
    \end{eqnarray}
Thus we can replace $F_1^t$ and $F_2^t$ with $F_1^k$ and $F_2^k$ in (\ref{splitting_scheme1}), based on which the adversary responses with the same answers of the queries with $(\z_1^t,\z_2^{t+1},\beta^t)$. In other words, this replacement does not influence the behavior of the algorithm and (\ref{splitting_scheme1}) produces the same sequence of $\{\x_1^k,\x_2^k\}$ with the following algorithm scheme, which performs
    \begin{eqnarray}
    \begin{aligned}\label{splitting_scheme2}
        &\mbox{Generate } \z_1^t\mbox{ in the same way with (\ref{splitting_scheme1})},\\
        &\x_1^{t+1}=\mbox{Prox}_{F_1^{k}/\beta^t}(\z_1^t),\quad F_1(\z_1^t)=F_1^k(\z_1^t),\\
        &\mbox{Generate } \z_2^{t+1} \mbox{ in the same way with (\ref{splitting_scheme1})},\\
        &\x_2^{t+1}=\mbox{Prox}_{F_2^{k}/\beta^t}(\z_2^{t+1}),\quad F_2(\z_2^{t+1})=F_2^k(\z_2^{t+1}),
    \end{aligned}
    \end{eqnarray}
at the $t$-th iteration. In scheme (\ref{splitting_scheme2}), we use $F_1^{k}$ and $F_2^{k}$, rather than $F_1^{t}$ and $F_2^{t}$.

We can prove that $F_i^k$ is convex and 1-Lipschitz continuous. $F^k(\x)\equiv F_1^k(\x)+F_2(\x)$ achieves the minimum at $\x^*=b\sum_{r=0}^{2k}\v_r$. If we let $b=\frac{1}{\sqrt{2k+1}}$, then $\|\x^*\|=1$. Due to the special form of $F^k(\x)$, we can prove $F^{k}(\hat\x^k)-F^{k}(\x^*)\geq \frac{1}{8(k+1)}$.

$\hfill\blacksquare$
\subsection{General ADMM Type Methods}
Now we use Proposition \ref{splitting theorem} to establish the lower complexity bound of ADMM type methods. Consider the following special case of problem (\ref{problem}):
\begin{equation}\label{problem_special}
\min_{\x_1,\x_2\in\mathbb X} F_1(\x_1)+F_2(\x_2),\quad s.t. \quad\x_1-\x_2=0.
\end{equation}
We consider the general ADMM type methods with the property of alternatingly minimizing the augmented Lagrangian function. Specifically, define the general ADMM type methods as
\begin{eqnarray}
\begin{aligned}\label{general_adm}
&\mbox{Generate } \lambda_2^t\mbox{ based on }\{\x_1^{1:t},\x_2^{1:t},\lambda_1^{1:t},\lambda_2^{1:t-1}\}\mbox{ and }\y_2^t\mbox{ based on }\{\x_1^{1:t},\x_2^{1:t}\},\\
&\x_1^{t+1}=\argmin_{\z} L(\x_1,\y_2^t,\lambda_2^t,\beta^t)=\mbox{Prox}_{F_1/\beta^t}\left(\y_2^{t}-\frac{\lambda_2^{t}}{\beta^t}\right),\\
&\mbox{Generate } \lambda_1^{t\hspace*{-0.03cm}+\hspace*{-0.03cm}1}\mbox{ based on }\{\hspace*{-0.05cm}\x_1^{1\hspace*{-0.02cm}:\hspace*{-0.02cm}t\hspace*{-0.03cm}+\hspace*{-0.03cm}1}\hspace*{-0.07cm},\hspace*{-0.05cm}\x_2^{1\hspace*{-0.02cm}:\hspace*{-0.02cm}t}\hspace*{-0.07cm},\hspace*{-0.05cm}\lambda_1^{1\hspace*{-0.02cm}:\hspace*{-0.02cm}t}\hspace*{-0.07cm},\hspace*{-0.05cm}\lambda_2^{1\hspace*{-0.02cm}:\hspace*{-0.02cm}t}\hspace*{-0.05cm}\}\mbox{ and }\y_1^{t\hspace*{-0.03cm}+\hspace*{-0.03cm}1}\mbox{ based on }\{\hspace*{-0.05cm}\x_1^{1\hspace*{-0.02cm}:\hspace*{-0.02cm}t\hspace*{-0.03cm}+\hspace*{-0.03cm}1}\hspace*{-0.07cm},\hspace*{-0.05cm}\x_2^{1\hspace*{-0.02cm}:\hspace*{-0.02cm}t}\hspace*{-0.07cm}\},\\
&\x_2^{t+1}=\argmin_{\x} L(\y_1^{t+1},\x_2,\lambda_1^{t+1},\beta^t)=\mbox{Prox}_{F_2/\beta^t}\left(\y_1^{t+1}-\frac{\lambda_1^{t+1}}{\beta^t}\right),
\end{aligned}
\end{eqnarray}
at the $t$-th iteration and $\beta^t$ can be any value. It can be checked that the traditional ADMM and ALADMM-NE (with $f_i=0$ and $\A_i=\I$) belong to this general scheme.

We can see that procedure (\ref{general_adm}) belongs to (\ref{splitting_scheme}) by letting $\z_1^t=\y_2^t-\frac{\lambda_2^{t}}{\beta^t}$ and $\z_2^{t+1}=\y_1^{t+1}+\frac{\lambda_1^{t+1}}{\beta^t}$. Letting $\hat\x_1^k=\sum_{i=1}^k\alpha_1^i\x_1^i$ and $\hat\x_2^k=\sum_{i=1}^k\alpha_2^i\x_2^i$, then from Proposition \ref{splitting theorem} we know that there exists convex and $L$-continuous $F_1$ and $F_2$ such that $F_1(\hat\x_2^k)+F_2(\hat\x_2^k)-F_1(\x^*)-F_2(\x^*)\geq \frac{LB}{8(k+1)}$. Since $F_1$ is $L$-continuous: $|F_1(\hat\x_2^k)-F_1(\hat\x_1^k)|\leq L\|\hat\x_2^k-\hat\x_1^k\|$, we can have $F_1(\hat\x_2^k)\leq F_1(\hat\x_1^k)+L\|\hat\x_2^k-\hat\x_1^k\|$ and
\begin{subequations} \begin{align}
\frac{LB}{8(k+1)}\leq&F_1(\hat\x_2^k)+F_2(\hat\x_2^k)-F_1(\x^*)-F_2(\x^*)\notag\\
\leq& L\|\hat\x_2^k-\hat\x_1^k\|+F_1(\hat\x_1^k)+F_2(\hat\x_2^k)-F_1(\x^*)-F_2(\x^*)\notag\\
\leq& L\|\hat\x_2^k-\hat\x_1^k\|+|F_1(\hat\x_1^k)+F_2(\hat\x_2^k)-F_1(\x_1^*)-F_2(\x_2^*)|\notag
\end{align} \end{subequations}
where $\x^*=\x_1^*=\x_2^*$. Thus we have the following lower complexity bound proposition for the general ADMM type methods for both the ergodic and nonergodic case, where the nonergodic bound can be obtained by letting $\alpha_1^i=\alpha_2^i=0, i=1,\cdots,k-1,$ and $\alpha_1^k=\alpha_2^k=1$.

\begin{proposition}\label{splitting theorem1}
For any algorithm belonging to the general splitting scheme (\ref{general_adm}), there exists convex and $L$-continuous functions $F_1$ and $F_2$ defined over $\mathbb X=\{\x\in \mathbb{R}^{6k+2}:\|\x\|\leq B\}$, such that
\begin{subequations} \begin{align}
L\|\hat\x_2^k-\hat\x_1^k\|+|F_1(\hat\x_1^k)+F_2(\hat\x_2^k)-F_1(\x_1^*)-F_2(\x_2^*)|\geq \frac{LB}{8(k+1)}.\notag
\end{align} \end{subequations}
where $\hat\x_1^k=\sum_{i=1}^k\alpha_1^i\x_1^i$ and $\hat\x_2^k=\sum_{i=1}^k\alpha_2^i\x_2^i$, $\forall\alpha_1^i$ and $\forall\alpha_2^i$, $i=1,\cdots,k$.
\end{proposition}

Since problem (\ref{problem_special}) is a special case of problem (\ref{problem}), we can have that $O(1/K)$ is the optimal convergence rate of the general ADMM type methods (\ref{general_adm}) for problem (\ref{problem}). There is no better ADMM type algorithm which converges faster than the $O(1/K)$ rate if it belongs to the framework of (\ref{general_adm}). Moreover, (\ref{general_adm}) is general enough for the separable problem (\ref{problem}) while still keeping the property of ADMM that alternately minimizes the augmented Lagrangian function. Thus our result is general enough. Since we can easily construct some algorithms (which may diverge) such that they can easily make one of $\|\A\x-\b\|$ and $|F(\x)-F(\x^*)|$ small but difficult to keep both small, this is why we use the summation in Proposition \ref{splitting theorem1}.

\section{Experiments on the Group Sparse Logistic Regression with Overlap}
In this section we test the performance of ALADMM-NE and ALADMM-NER on the Group Sparse Logistic Regression with Overlap. This problem can be deemed as a combination of the Group Sparse Logistic Regression \cite{Meier-2008-GroupLogit} and the Group LASSO with Overlap \cite{Jacob-2009-GroupLASSO}. Its mathematical model is as follows:
\begin{subequations} \begin{align}
\hspace*{-0.1cm}\min_{\w,\b} \frac{1}{s}\sum_{i=1}^s \mbox{log}(1+\mbox{exp}(-\y_i(\w^T\x_i+\b)))+\nu\sum_{j=1}^t \|\S_j\w\|,\notag
\end{align} \end{subequations}
where $\x_i$ and $\y_i$ are the training samples and labels. $\w$ and $\b$ are the parameters for the classifier. $s$ is the sample size and $t$ is the group size. $\S_j,j=1,\cdots,t$ are the selection matrices with only one 1 at each row and 0 for the rest entries. We consider the case that the groups of entries may overlap each other. We can transform the problem to a linearly constrained one by introducing $\overline \S_j=(\S_j;\mathbf 0)$, $\overline\S=\left(
  \begin{array}{c}
     \overline\S_1\\
     \vdots\\
     \overline\S_t
  \end{array}
\right)$, $\overline\w=\left(
  \begin{array}{c}
     \w\\
     \b
  \end{array}
\right)$, $\overline\x_i=\left(
  \begin{array}{c}
     \x\\
     1
  \end{array}
\right)$, $\z_j= \overline\S_j\overline\w$ and $\z=\left(
  \begin{array}{c}
     \z_1\\
     \vdots\\
     \z_t
  \end{array}
\right)$:
\begin{equation}\label{logit_model}
\min_{\overline\w,\z} \frac{1}{s}\sum_{i=1}^s \mbox{log}(1+\mbox{exp}(-\y_i(\overline\w^T\overline\x_i)))+\nu\sum_{j=1}^t \|\z_j\|,\quad s.t.\quad \z= \overline\S\overline\w.
\end{equation}

\begin{figure}
\centering
\begin{tabular}{@{\extracolsep{0.001em}}c@{\extracolsep{0.001em}}c}
\includegraphics[width=0.5\textwidth,keepaspectratio]{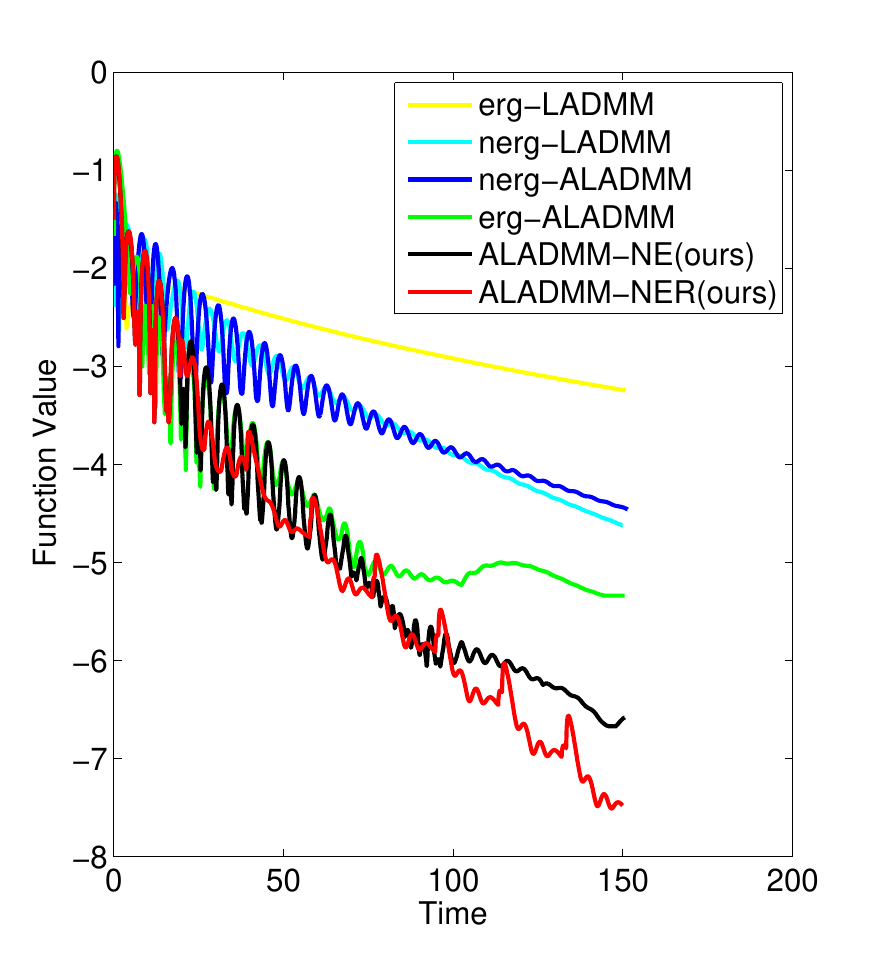}
&\includegraphics[width=0.5\textwidth,keepaspectratio]{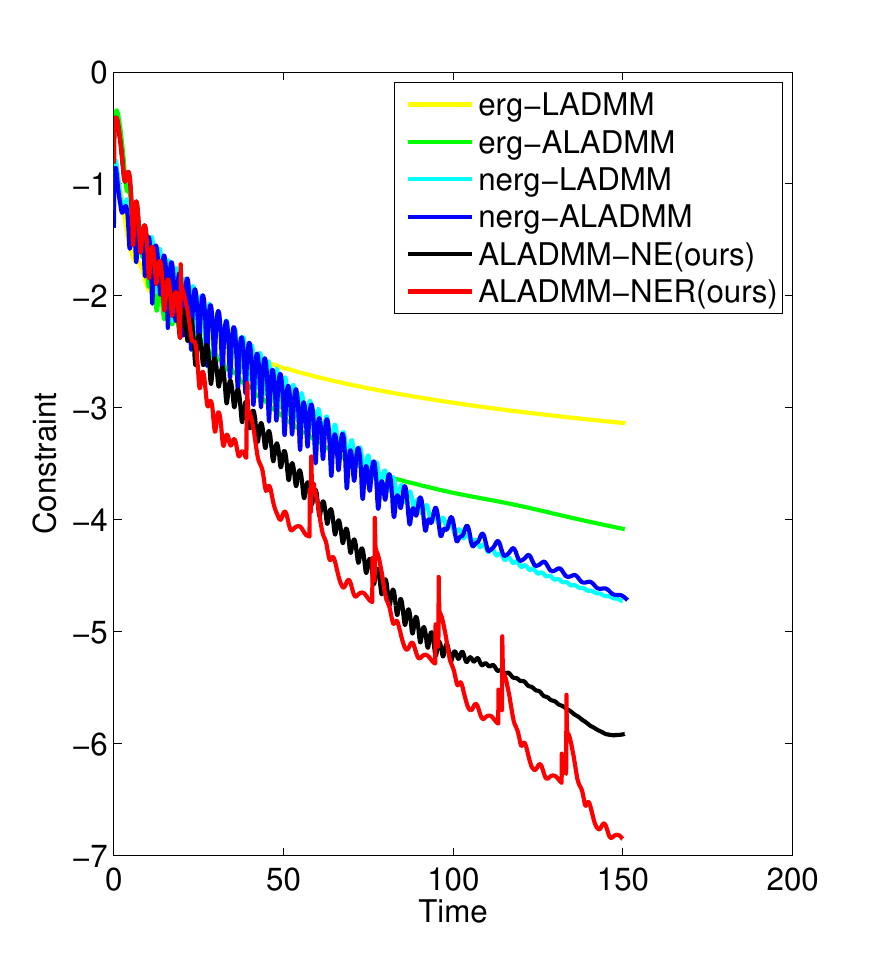}\\
(a) $\mbox{log}_{10}\left|F(\overline\w^k,\z^k)-F(\overline\w^*,\z^*)\right|$ & (b) $\mbox{log}_{10}\left\|\z^k-\overline\S\overline\w^k\right\|$\\
\includegraphics[width=0.5\textwidth,keepaspectratio]{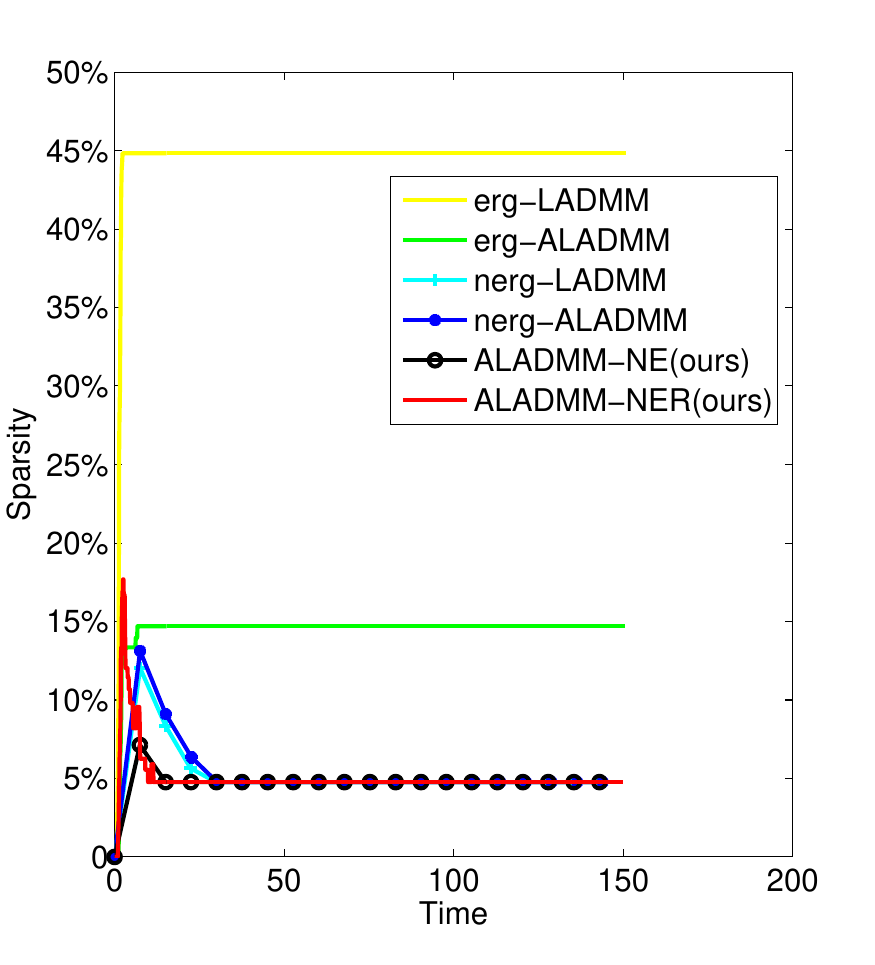}
&\includegraphics[width=0.5\textwidth,keepaspectratio]{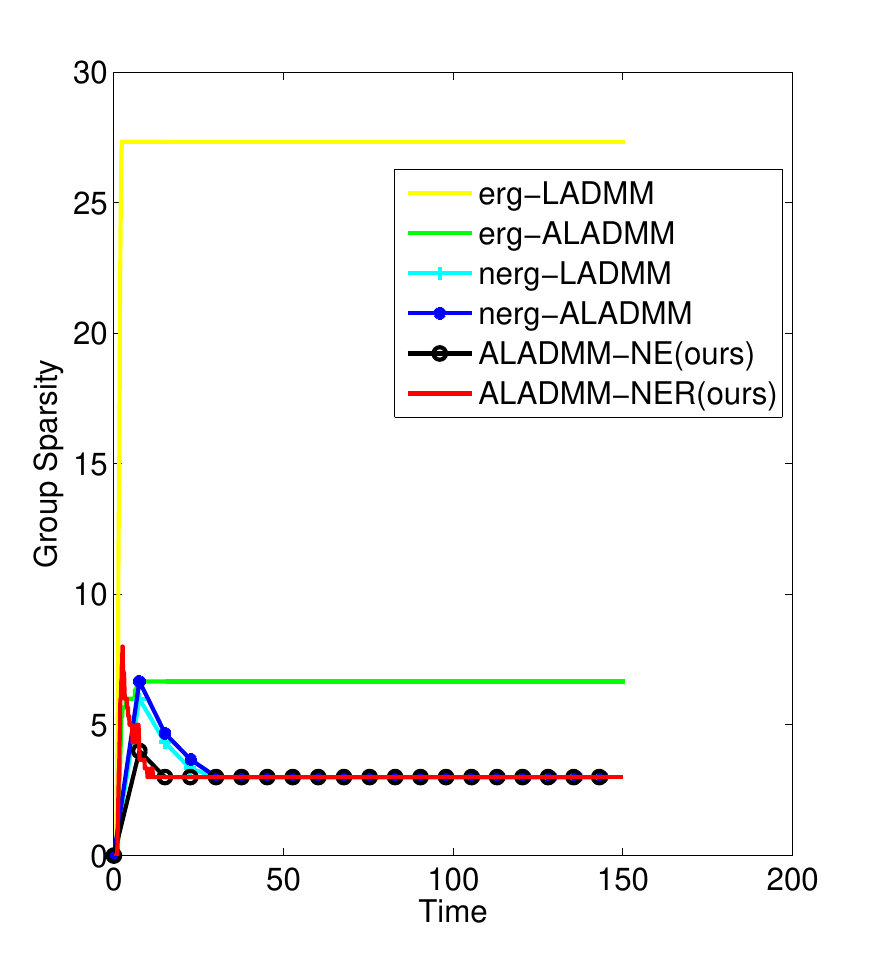}\\
(c) Sparsity & (d) Group Sparsity
\end{tabular}
\caption{Compare ALADMM-NE and ALADMM-NER with LADMM and ALADMM on the Group Sparse Logistic Regression problem. We present the function value, constraint error, sparsity (percent of selected Genes) and Group sparsity (number of non-empty groups).}\label{fig_lr}
\end{figure}
We carry out the experiment on the breast cancer gene expression data set \cite{Van-2002-Breast}. 3510 genes in 295 breast cancer tumors are considered in our experiment, which appear in 637 gene groups. Gene selection is a key purpose in this problem. The group sparsity regularization helps to decide which groups of Genes play a central role in the cancer prediction. Thus the group sparsity is strongly required.

We compare ALADMM-NE and ALADMM-NER with LADMM and the accelerated LADMM (ALADMM) \cite{ouyang-2015-fastLADM}. We set the initializer at 0 and run all the methods for 2000 iterations. We set $\tau=0.8$ for ALADMM-NE and ALADMM-NER and $\epsilon=0.02$ for ALADMM-NER. For ALADMM, we set the parameters following the assumptions in Theorem 2.6 of \cite{ouyang-2015-fastLADM}. We set $\beta=0.3$ for LADMM, $\beta=0.06$ for erg-ALADMM, $\beta=0.4$ for nerg-ALADMM, $\beta=0.08$ for ALADMM-NE and ALADMM-NER for the best performance of each algorithm, respectively, where erg-ALADMM (erg-LADMM) means that we use the ergodic solution $\x^K$ for ALADMM ($\sum_{k=1}^K \x^k/K$ for LADMM) and nerg-ALADMM (nerg-LADMM) means that we use the nonergodic solution $\z^K$ for ALADMM ($\x^K$ for LADMM).Ouyang et al. \cite{ouyang-2015-fastLADM} proposed a backtracking scheme to estimate $\|\overline\S\|_2$ and the Lipschitz constant $L$. Since $\|\overline\S\|_2$ and $L$ can be exactly computed in our problem, we do not use the backtracking scheme for simplicity.

Figure \ref{fig_lr} draws the plots of the objective function value, the constraint error, the sparsity and the group sparsity vs. time. We run LADMM for 100000 iterations and use its nonergodic output as the optimal $(\overline\w^*,\z^*)$, which is used to plot $\left|F(\overline\w^k,\z^k)-F(\overline\w^*,\z^*)\right|$. We can see that both erg-LADMM and erg-ALADMM have a less favorable sparsity and group sparsity than their nonergodic counterparts, this verifies that the nonergodic measurement is required. However, Nerg-ALADMM decreases the objective function slower than erg-ALADMM. In some practical applications, ADMM can perform better than the theoretical bound. Thus it is not strange that nerg-LADMM converges faster than erg-LADMM. As comparison, ALADNM-NE and ALADMM-NER not only run faster than the compared methods but also have the sparsity and group sparsity as well as nerg-LADMM and nerg-ALADMM. In ADMM type methods, the monotonicity of the objective function and the constraint error cannot be guaranteed in theory. This leads to the oscillation in Figure \ref{fig_lr}.

\section{Conclusions}
In this paper, we modify the accelerated ADMM proposed in \cite{ouyang-2015-fastLADM} and give an $O(1/K)$ nonergodic analysis in the sense of $|F(\x^K)-F(\x^*)|\leq O(1/K)$ and $\|\A\x^K-\b\|\leq O(1/K)$, where the nonergodic result has a more favorable sparseness and low-rankness than the ergodic one. This is the first $O(1/K)$ nonergodic convergent ADMM type method and surpasses the $o(1/\sqrt{K})$ nonergodic rate of the traditional ADMM. Moreover, we show that the lower complexity bound of ADMM type methods is $O(1/K)$ when each $F_i$ is nonsmooth and non-strongly convex, which means that our method is optimal.

\section{Acknowledgement}

Zhouchen Lin is supported by National Basic Research Program of China (973 Program) (grant no. 2015CB352502), National Natural Science Foundation (NSF) of China (grant nos. 61625301 and 61731018), Qualcomm and Microsoft Research Asia.

\small
\bibliographystyle{unsrt}
\bibliography{arxiv}
\end{document}